\documentclass[11pt]{amsart}

\usepackage[utf8x]{inputenc} 	     
\usepackage{lmodern}                
\usepackage[T1]{fontenc}            
\usepackage[english]{babel}         
\usepackage{ucs}             	    
\usepackage{hyperref}               

\usepackage{geometry}
\usepackage{amsmath}
\usepackage{amssymb}
\usepackage{bm}
\usepackage{graphicx}
\usepackage{amsthm}
\usepackage{color}
\usepackage{tabu}
\usepackage{booktabs}
\usepackage{caption}

\newtheorem{theorem}{Theorem}

\newtheorem{rem}{Remark}

\newcommand{\RR}{\mathbb{R}}
\newcommand{\SSS}{\mathbb{S}}
\newcommand{\Q}{\mathcal Q}
\newcommand{\I}{\mathcal I}
\newcommand{\ve}{\varepsilon}
\newcommand{\diverg}{\nabla \cdot}

\newcommand*\diff{\mathop{}\!\mathrm{d}}

\title[Granular Materials]{Recent development in kinetic theory of granular materials: analysis and numerical methods}
\author{Jose Antonio Carrillo}
\address{\flushleft
  \underline{Jose Antonio Carrillo}\\[.3em]
  Department of Mathematics, Imperial College London\\[.2em]
  SW7 2AZ London, UK\\[.2em]
  Email: \texttt{carrillo@imperial.ac.uk}
}
\author{Jingwei Hu}
\address{\flushleft  
  \underline{Jingwei Hu}\\[.3em]
  Department of Mathematics, Purdue University\\[.2em]
  West Lafayette, IN 47907, USA\\[.2em]
  Email: \texttt{jingweihu@purdue.edu}
}
\author{Zheng Ma}
\address{\flushleft  
  \underline{Zheng Ma}\\[.3em]
  Department of Mathematics, Purdue University\\[.2em]
  West Lafayette, IN 47907, USA\\[.2em]
  Email: \texttt{ma531@purdue.edu}
}
\author{Thomas Rey}
\address{\flushleft
  \underline{Thomas Rey}\\[.3em]
  Univ. Lille, CNRS, UMR 8524, Inria – Laboratoire Paul Painlevé\\[.2em]
  F-59000 Lille, France\\[.2em]
  Email: \texttt{thomas.rey@univ-lille.fr}}

\thanks{JAC acknowledges support by the Engineering and Physical Sciences Research Council (EPSRC) under grant no. EP/P031587/1 and by the National Science Foundation (NSF) under grant no. RNMS11-07444 (KI-Net). JH was partially funded by NSF grant DMS-1620250 and NSF CAREER grant DMS-1654152. TR was partially funded by Labex CEMPI (ANR-11-LABX-0007-01) and ANR Project MoHyCon (ANR-17-CE40-0027-01). }

\begin{document}

\begin{abstract}
Over the past decades, kinetic description of granular materials has received a lot of attention in mathematical community and applied fields such as physics and engineering. This article aims to review recent mathematical results in kinetic granular materials, especially for those which arose since the last review \cite{Villani:2006} by Villani on the same subject. We will discuss both theoretical and numerical developments. We will finally showcase some important open problems and conjectures by means of numerical experiments based on spectral methods.\\[3em]
\end{abstract}

\maketitle
\tableofcontents

\section{The Boltzmann equation for granular gases}

Granular gases have been initially introduced to describe the nonequilibrium behavior of materials composed of a large number of non necessarily microscopic particles, such as grains or sand. These particles form a gas, interacting \textit{via} energy dissipating inelastic collisions. Although a statistical mechanics description of particle systems through inelastic collisions faces basic derivation problems as the inelastic collapse \cite{McNamaraYoung:1993PRE}, i.e. infinite many collisions in finite time, the kinetic description of rapid granular flows \cite{JR85,JR852,Goldhirsch:2003} has been able to compute transport coefficients for hydrodynamic descriptions used successfully in situations that are a long way from their supposed limits of validity, to describe, for instance, shock waves in granular gases \cite{RBSS02,BMSS02}, clustering \cite{HM03,BSSP04,CS}, and the Faraday instability for vibrating thin granular layers \cite{Faraday:1831,MeloUmbanhowarSwinney:1995,MeloUmbanhowarSwinney:1996,TsimringAranson:1997,UmbanhowarMeloSwinney:1998,BMSS02,Bougieetal04,carrillo:2008granular}. A large amount of practical systems can be described as a granular gas, such as for example spaceship reentry in a dusty atmosphere (Mars for instance), planetary rings \cite{Kawai:1990,Araki:1986} and sorting behavior in vibrating layers of mixtures. A lot of other examples can be found in the thesis manuscript \cite{Daerr:2000}, and in the seminal book of Brilliantov and P\"oschel \cite{brilliantov:2004}. 

Usually, a granular gas is composed of $10^6$ to $10^{16}$ particles. The study of such a system will then be impossible with a direct approach, and we shall adopt a kinetic point of view, studying the behavior of a one-particle distribution function $f$, depending on time $t\geq 0$, space $x \in \Omega \subset \mathbb R^{d_x}$ and velocity $v \in \mathbb R^{d}$, for $d_x \leq d \in \{1,2,3\}$. The statistical mechanics description of the system has been then admitted in the physical community as the tool to connect the microscopic description to macroscopic system of balance laws in rapid granular flows \cite{JR85,JR852,Goldhirsch:2003,GS95,BDS99} as in the classical rarefied gases \cite{CIP:94}. In this first section, we shall review some basics on the inelastic Boltzmann equation, and present the mathematical state of the art since the previous review paper on the subject \cite{Villani:2006}.

\smallskip	  
	  \paragraph{\textbf{Microscopic dynamics.}}
	   
\begin{figure}[ht]
			\begin{center}
				\includegraphics{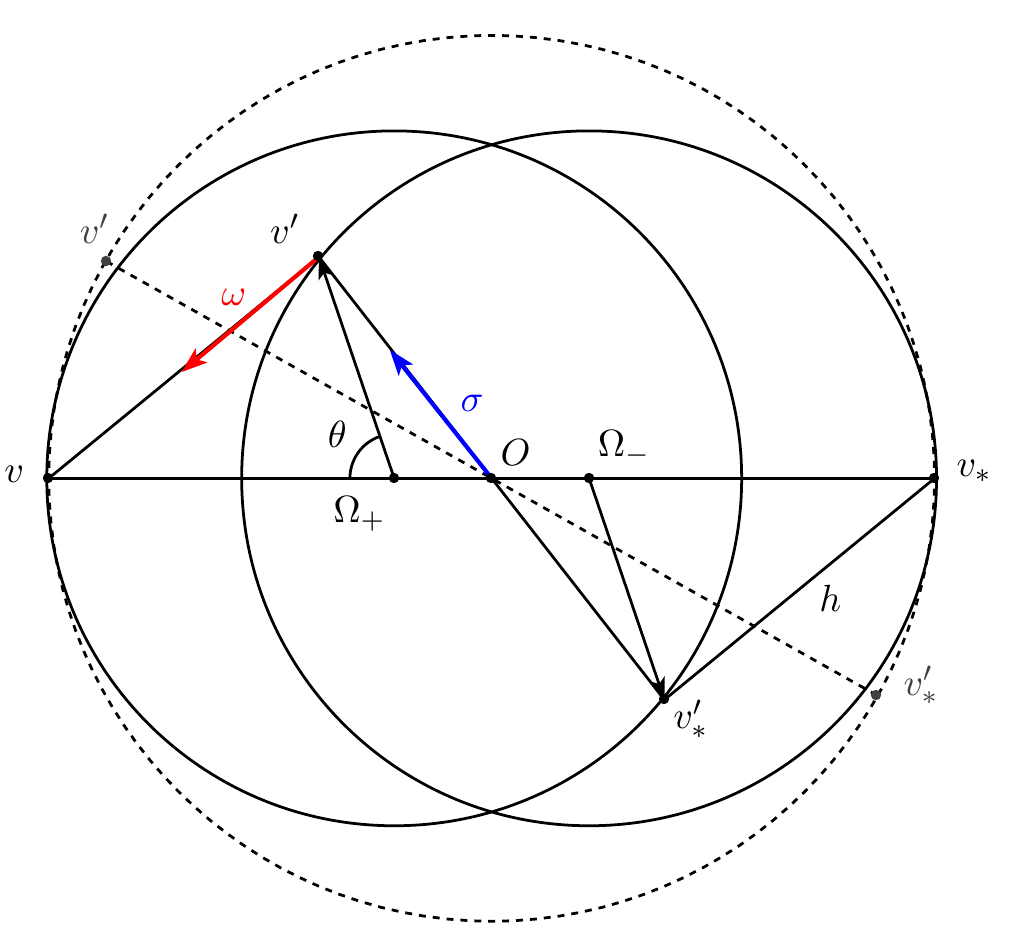}
				\caption{Geometry of the inelastic collision in the phase space (dashed lines represent the elastic case).}
				\label{figCollSphereInel}
			\end{center}
		\end{figure}

	   The microscopic dynamics can be summarized with the following hypotheses:
	    \begin{enumerate}
	    
	       \item \label{hypColBin}
	        The particles interact via  \emph{binary} collisions. More precisely, the gas is rarefied enough so that collisions between 3 or more particles can be neglected.
	        
	      \item \label{hypColLocal}
	        These binary collisions are localized in space and time. In particular, all the particles are considered as point particles, even if they describe macroscopic objects. 
	        
	      \item \label{hypColConserv}
	        Collisions preserve mass and momentum, but dissipate a fraction $1-e$ of the kinetic energy in the impact direction, where the inelasticity parameter $e \in [0,1]$ is called \emph{restitution coefficient}:
	      \begin{equation} 
					\label{eqMicroConsInel}
					\left\{ \begin{aligned}
					& v' \,+\, v_*' \,=\, v \,+\, v_*, \\
					& |v'|^2 \,+\, |v_*'|^2 \,-\, |v|^2 \,-\, |v_*|^2  \,=\, - \frac{1-e^2}{2} \,| (v-v_*) \cdot \omega |^2 \,\leq\, 0,
					\end{aligned} \right.
			  \end{equation}
		      with  $\omega \in \SSS^{d-1}$ being the impact direction.
            Using these conservation, one has the following two possible parametrizations (see also Fig.~\ref{figCollSphereInel}) of the post-collisional velocities, as a function of the pre-collisional ones:
			  \begin{itemize}
				  \item The $\omega$-representation or reflection map, given for $\omega \in \SSS^{d-1}$ by
						\begin{align}
						    \label{eq:omegaRep}
						  v' &= v - \frac{1 + e}{2} \left( (v - v_*) \cdot \omega \right) \omega, \nonumber \\
				    	v_*' &= v_* + \frac{1 + e}{2} \left( (v - v_*) \cdot \omega \right) \omega.
						\end{align}
    	         \item The $\sigma$-representation or swapping map, given for $\sigma \in \SSS^{d-1}$ by
		        \begin{align} 
		            \label{eq:sigmaRep}
					v' &= \frac{v + v_*}{2} + \frac{1 - e}{4} (v-v_*) + \frac{1+e}{4} |v - v_*| \sigma, \nonumber \\
				    v_*' &= \frac{v + v_*}{2} - \frac{1 - e}{4} (v-v_*) - \frac{1+e}{4} |v - v_*| \sigma.
				\end{align}
				\label{hypColInel}
		      \end{itemize}
	    \end{enumerate}

\begin{figure}[ht]
			\begin{center}
				\includegraphics[scale=1.1]{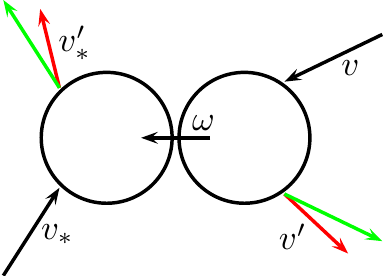}
				\caption{Geometry of the inelastic collision in the physical space (green is elastic, red is inelastic).}
				\label{fig:collisionElInel}
			\end{center}
		\end{figure}

        \begin{rem}
            Taking $e = 1$ in both \eqref{eq:omegaRep} and \eqref{eq:sigmaRep} yields the classical energy-conservative elastic collision dynamics, as illustrated in Fig.~\ref{fig:collisionElInel}.
        \end{rem}
        
		The geometry of collisions is more complex than the classical elastic one.
		Indeed, fixing $v, v_* \in \mathbb{R}^d$, denote by  
		\[ 
		    \Omega_{\pm} :=\frac{v+v_*}{2} \pm \frac{1-e}{2} (v_* - v), \qquad O := \frac{v+v_*}{2} = \frac{v'+v_*'}{2}.
		\]
		Then if $u := v - v_*$ is the  \emph{relative velocity}, one has 
		\[
		    |\Omega_- - v' | = | \Omega_+ - v_*'| = \frac{1+e}{2} |u|, 
		\] 
		namely $v' \in \mathcal{S} \left(\Omega_-, |u|(1+e)/2  \right)$ and $v_*' \in \mathcal{S} \left( \Omega_+, |u|(1+e)/2 \right) $, where $\mathcal{S}(x, r)$ is the sphere centered in $x$ and of radius $r$ (see also Fig.~\ref{figCollSphereInel}).

		\smallskip
	    \paragraph{\textbf{Restitution coefficient.}}
	    The physics literature is quite divided on the question of whether the restitution coefficient $e$ should be a constant or not \cite{brilliantov:2004}. 
	    Although most of the early mathematical results on the topic consider a constant $e$ \cite{Villani:2006}, it seems that this case is only realistic in dimension 1 of velocity (the so-called ``collisional cannon'' described in the chapter 4 of \cite{brilliantov:2004} is a famous counter-example).
	    The true realistic case considers that  $e$ depends on the relative velocity  $|v-v_*|$  of the colliding particles. 
	    Even more precisely, it must be close to the elastic case $1$ for small relative velocities (namely no dissipation, elastic case), and decay towards $0$ when this relative velocity is large. The first mathematical result on this direction can be found in \cite{Toscani:2000}, where
	    \begin{equation}
	      \label{defRestiCoeffToscani}
	      e(|v-v_*|) = \dfrac{1}{1+c\,|v-v_*|^\gamma},
	    \end{equation} 
	    for a nonnegative constant $c$ characterizing the inelasticity strength ($c=0$ being elastic), and $\gamma \in \RR$. 
	    
	    Another important case is the so-called  \emph{viscoelastic} hard spheres one, thoroughly studied mathematically in a series of papers \cite{alonso:2009,alonso2013uniqueness,bisi2011uniqueness,alonso2014boltzmann}, where $e$ is given by the implicit relation
	    \begin{equation} 
	      \label{defRestiCoeffViscoEl}
	      e(|v-v_*|) + a |v-v_*|^{1/5} e(|v-v_*|)^{3/5} = 1, 
	    \end{equation}
	    for $a>0$.
	    More details on the derivation of this expression can be found in \cite[Chapter 3]{brilliantov:2004}.
	    
	    Another quite rigorous study has been made in \cite{poschel:2005transient}, with a threshold-dependent restitution coefficient:
	    \[ e(r) = \left \{ \begin{aligned} 1 && \text{ if } r < r_*, \\ \bar e && \text{ if } r \geq r_*, \end{aligned} \right. \]
	    $\bar e < 1$ and $r_* > 0$ being fixed. 
	    
	    Finally, the case $e = 0$ describes \emph{sticky collisions}: the normal component of the kinetic energy being completely dissipated during impact, the particles stick and travel together in the tangent direction after impact. A derivation of the model from the microscopic dynamics on the line can be found in \cite{BrenierGrenier:1998,ERykovSinai:1996}.
	    
	  	\begin{rem}
			This model is meaningful even in dimension $1$, which is not the case for elastic collisions. Indeed, such monodimensional collisions are only
			\[\{v', v_*' \} = \{ v, v_* \},\]
			meaning that the particle velocities are either swapped or preserved. The particles being indistinguishable, nothing happens\footnote{Because of that, the elastic collision operator is simply equal to $0$ for a one-dimensional velocity space, the Boltzmann equation reducing only to the free transport equation.}. 
			In the $1$d inelastic case, the collisions are given using \eqref{eq:sigmaRep} by 
			\[ \{v', v_*' \} = \{ v, v_* \}  \quad \text{ or } \quad \left\{ \frac{v+v_*}{2} \pm \frac{e}{2} (v - v_*) \right\} \] 
			depending on the value of  $\sigma \in \{ \pm1 \}$.
		\end{rem}

		\paragraph{\textbf{The granular gases operator: Weak form.}	}	
			Using the microscopic hypotheses (\ref{hypColBin}--\ref{hypColLocal}--\ref{hypColInel}), one can derive the granular gases collision operator $\Q_\mathcal{I}$, by following the usual elastic procedure (see \cite{Villani:2006} for more details). Its \emph{weak form} in the $\sigma$-representation is given by
			\begin{equation} 
			  \label{defQGranWeak}
				\int_{\RR^d} \Q_\mathcal{I}(f,f)(v) \,\psi(v)\,dv = \frac{1}{2 }\int_{\mathbb{R}^d \times \mathbb{R}^d \times \mathbb{S}^{d-1}} f_{*} \, f \, \left(\psi' + \psi_*' - \psi - \psi_* \right) B(|v-v_*|,\cos \theta, E(f)) \, d\sigma \, dv \, dv_*,
			\end{equation}	
where  the collision kernel is typically of the form $B(|u|,\cos \theta, E(f))=\Phi(|u|) b(\cos \theta, E(f))$, and $E(f)$ is the kinetic energy of $f$, namely its second moment in velocity, the postcollisional velocities are computed by \eqref{eq:sigmaRep}, and $\theta$ is the angle between $\sigma$ and $u$.
We shall assume in all the following of this section that the collision kernel is of \emph{generalized hard sphere} type, namely 
			 \begin{equation}
		        \label{defHSkernelAnomalous}
		        B(|u|,\cos \theta, E))=\Phi(|u|)\widehat{b}(\cos \theta, E ) = |u|^\lambda \, b(\cos \theta) \, E^\gamma,
		    \end{equation}
		    where $\lambda \in [0,1]$ ($\lambda = 0$ being the simplified \emph{Maxwellian pseudo-molecules case} and $\lambda = 1$ the more relevant \emph{hard sphere} case),  $\gamma \in \RR$ and the angular cross section $b$ verifies
		    \begin{equation} 
		        \label{def:AngularCrossSection}
		        0 < \beta_1 \leq b(x) \leq \beta_2 < \infty, \quad \forall x \in [-1, 1]. 
		    \end{equation}
		    
		    \begin{rem}
		        Note that we assumed that the collision kernel $B$ in \eqref{defQGranWeak} depends on the relative velocity, the angle of collision, and on $E(f)$. 
		        These former dependencies are quite classical, but the latter is not. Nevertheless, it makes a lot of sense physically speaking, as one can see in \cite{rey:2011}.
		    \end{rem}
The weak form in the $\omega$-representation can be written analogously as
\begin{equation} 
			  \label{defQGranWeakomega}
				\int_{\RR^d} \Q_\mathcal{I}(f,f)(v) \,\psi(v)\,dv = \frac{1}{2 }\int_{\mathbb{R}^d \times \mathbb{R}^d \times \mathbb{S}^{d-1}} f_{*} \, f \, \left(\psi' + \psi_*' - \psi - \psi_* \right) \widetilde B(|u|, \cos \theta, E(f) ) \, d\omega \, dv \, dv_*,
			\end{equation}	
where the postcollisional velocities are computed by \eqref{eq:omegaRep}, $\theta$ is the angle between $\omega$ and $u$, and 
\begin{equation*}
    \widetilde B(|u|,\cos \theta, E)=|u|^\lambda \widetilde b(\cos \theta) E^\gamma
\end{equation*} 
with $\widetilde b(t) =  3|t| b(1-2t^2)$ for $-1 \le t \le 1$ by the change of variables between the $\sigma$- and the $\omega$-representation, see \cite{CCC} for details.

        \smallskip
		\paragraph{\textbf{The granular gases operator: Strong form.}}	
		    
		    Deriving a strong form of $\Q_\mathcal{I}$ with the reflection map in the $\omega$-representation is a matter of a change of variables. However,
		    deriving a strong form of $\Q_\mathcal{I}$ is not as easy as in the elastic case in the $\sigma$-representation since the collisional transform $(v,v_*, \sigma) \to (v',v_*',\sigma)$ is not an involution and we have to go through the $\omega$-representation, see \cite{CCC} for details. 
		    
		    More precisely, given the restitution coefficient $e = e(|u|)$ depending on the relative velocity of the particles $u=v-v_*$, we assume the collisional transform's Jacobian for \eqref{eq:omegaRep} is $J(|u|,\cos \theta)\neq 0$ for all $z$. Notice $J=e$ in the constant restitution case. It is in general  a complicated expression of the relative speed $r=|u|$ and $s=\cos \theta$ involving $e$ and its derivative.
Then, the precollisional velocities read as
\begin{align}
						    \label{eq:omegaReppost}
						  'v &= v - \frac{1 + e}{2e} \left( (v - v_*) \cdot \omega \right) \omega, \nonumber \\
				    	'v_* &= v_* + \frac{1 + e}{2e} \left( (v - v_*) \cdot \omega \right) \omega.
\end{align}
The final strong from of the operator is $\Q_\mathcal{I}(f, f)(v)=\Q_\mathcal{I}^+(f,f)(v) - f(v) \, L(f)(v)$ with the loss part of the operator given by
$$
L(f)(v)= \int_{\mathbb{R}^d \times \mathbb{S}^{d-1}}  \widetilde B(|v-v_*|, \cos \theta, E(f) ) f_* \, d \omega \, dv_*
$$
and the gain part of the operator in strong form written as
\begin{align}
\Q_\mathcal{I}^+(f, f)(v) & = \int_{\mathbb{R}^d \times \mathbb{S}^{d-1}} \widetilde \Phi_e^+(|u|,\cos \theta) \widetilde b_e^+
(\cos \theta) \frac{E^\gamma}{J(|u|,\cos \theta)} \, \,'f \,'f_*\, d \omega \, dv_*, \label{defQGran} 
\end{align}
with $\widetilde \Phi_e^+(r,s)$ and $\widetilde b_e^+(s)$ given by
\begin{equation}\label{betildef}
\widetilde b_e^+(s) = \widetilde b\left(\frac{s}{\sqrt{e^2
+(1-e^2)s^2 }}\right) \ ,
\end{equation}
and 
\begin{equation}\label{phietildef}
\widetilde \Phi_e^+(r,s) = \Phi\left(\frac{r}{e}
\sqrt{e^2 +(1-e^2)s^2}\right)= \left(\frac{r}{e}
\sqrt{e^2 +(1-e^2)s^2}\right)^\lambda  .
\end{equation}
We can derive now the following \emph{strong form} of the collision operator also in the $\sigma$-representation by just changing variable in the operator from $\omega$ to $\sigma$, see \cite{CCC}, to find the expressions of the loss and the gain terms in the $\sigma$-representation:
$$
L(f)(v)= \int_{\mathbb{R}^d \times \mathbb{S}^{d-1}} B(|v-v_*|, \cos \theta, E(f) ) f_* \, d \sigma \, dv_*
$$
and 
\begin{align}
\Q_\mathcal{I}^+(f, f)(v) & = \int_{\mathbb{R}^d \times \mathbb{S}^{d-1}} \Phi_e^+(|u|,\cos \theta)  b_e^+(\cos \theta) \frac{E^\gamma}{J(|u|,\cos \theta)} \, \,'f \,'f_*\, d \sigma \, dv_*, \label{defQGransigma} 
\end{align}
with $\Phi_e^+(r,s)$ and $b_e^+(s)$ given by
\begin{equation}\label{betildef_sigma}
b_e^+(s) = b\left(\frac{(1+e^2)s - (1-e^2)}{(1+e^2) -
(1-e^2)s}\right) \frac{\sqrt{2}}{\sqrt{{(1+e^2) -
(1-e^2)s}}} \ ,
\end{equation}
and 
\begin{equation}\label{phietildef_sigma}
\Phi_e^+(r,s)=\Phi\left(\frac{r}{\sqrt{2}e} \sqrt{(1+e^2) -
(1-e^2)s}\right)=\left(\frac{r}{\sqrt{2}e} \sqrt{(1+e^2) -
(1-e^2)s}\right)^\lambda  .
\end{equation}
In these expressions, the precollisional velocities are given in the $\sigma$-representation by
 \begin{align} 
		            \label{eq:sigmaReppre}
					'v &= \frac{v + v_*}{2} + \frac{1 - e}{4e} (v-v_*) + \frac{1+e}{4e} |v - v_*| \sigma, \nonumber \\
				    'v_* &= \frac{v + v_*}{2} - \frac{1 - e}{4e} (v-v_*) - \frac{1+e}{4e} |v - v_*| \sigma.
				\end{align}

			The granular gases collision operator has then the same structure of the elastic Boltzmann operator under Grad's cutoff assumption, namely the difference between an \emph{inelastic} gain term
			$\Q_\mathcal{I}^+(f,f)$ 
			and a loss operator
			$f \, L(f)$, 
			which depends only on the chosen collision kernel, but not on the inelasticity.
			
			We shall call the following \emph{granular gases equation}, or the inelastic Boltzmann equation:
			\begin{equation}
	            \label{eq:GranularGases}
	            \frac{\partial f}{\partial t} + v \cdot \nabla_x f = \Q_\mathcal{I}(f,f).
	        \end{equation}
		    We shall denote its first fluid moments (resp. density, mean momentum, and kinetic energy) by 
		    \[
		        \left ( \rho(f), \, \rho (f) \bm{u} (f), \,E(f) \right ) := \int_{\RR^{d}} \left( 1, v, |v|^2/2 \right ) f(v) \, dv.
		    \]

		\begin{rem}
			There is another popular approach to describe granular gases, which uses an \emph{Enskog}-type collision operator. It is more relevant physically because it allows to keep the particles' radii $\delta$ positive, hence delocalizing the collision\footnote{Note that using a BBGKY approach \cite{gallagher2014newton} to derive \eqref{defQGran} is not expected to succeed, because among other problems the macroscopic size of the particles composing a granular gas is incompatible with the Boltzmann-Grad scaling assumption.}.
			The microscopic hypothesis \ref{hypColLocal} is in particular not valid anymore. The strong form of the collision operator in the constant restitution coefficient case is given by
			\begin{equation}
			    \label{def:EnskogOp}
			  \Q_\mathcal{E}(f, f)(x,v) = \delta^{d-1} \, \int_{\mathbb{R}^d \times \mathbb{S}^{d-1}} \left( \widetilde \Phi_e^+(|u|,\cos \theta) \widetilde b_e^+
(\cos \theta) \frac{G(\rho_+)}{e} \,'f_+ \,'f_* - G(\rho_-) f_-\, f_* \right) \, d \omega \, dv_*,
			\end{equation}
			where $\rho$ is the local density of $f$,  $\pm$ denotes for a given function  $g = g(x)$ the shorthand notation 
			\[ g_\pm(x) := g(x \pm \delta \, \omega), \]
			and $G$ is the \emph{local} collision rate (also known as the \emph{correlation} rate, see  \cite{Villani:2006}). 
			The global existence of renormalized solutions for the full granular gases equation \eqref{eq:GranularGases} with the collision operator \eqref{def:EnskogOp} has been established for both elastic and inelastic collisions in \cite{Esteban:1991}.
			Existence and $L^1(dx\, dv)$ stability of such solutions has been proved in  \cite{Wu:2010}, for close to vacuum initial datum. 
	  \end{rem}

	\subsection{Cauchy theory of the granular gases equation}
		
		\smallskip
		\paragraph{\textbf{The space homogeneous setting.}}
		
	    Most of the rigorous mathematical results concerning the granular gases equation are obtained in the space homogeneous setting, where $f = f(t,v)$ is the solution to
          	\begin{equation} 
          	    \label{eq:GranularGasesHomogeneous}
				\left\{ \begin{aligned} 
				  & \partial_t f = \frac{\Q_\mathcal{I}(f, f)}\ve, \\
				  & f(0, v) = f_{in}(v),
				\end{aligned}\right. 
			\end{equation}
	    for a given scaling parameter $\ve >0$.

        The first existence results for solution to \eqref{eq:GranularGasesHomogeneous} can be found in \cite{bobylev:2000, bobylev:2001,BCT1,BCT2,BC,ReviewMetrics}. 
        These works deal with the generalized Maxwellian pseudo-molecule kernel \eqref{defHSkernelAnomalous} $\lambda = 0$, $b \equiv 1$ and $\gamma = 1/2$, with a velocity dependent restitution coefficient $e = e(|v-v_*|)$.
        Such a model allows to use some Fourier techniques to deal with the collision operator, altogether with the correct large time behavior for the kinetic energy, the so-called Haff's cooling Law \eqref{eq:HaffLaw}, and the correct hydrodynamic limit \eqref{sys:compressibleEulerGranular}. 
        The main result of these works is the global well-posedness of the solutions to \eqref{eq:GranularGasesHomogeneous} in $L_2^1(\RR^3)$, the convergence towards equilibrium, and the contraction in different metrics for the equation.
        The proof relies in the careful study of the self-similar solutions to \eqref{eq:GranularGasesHomogeneous}.
        Some extensions of these results, using the same collision kernel, can be found in the works \cite{bobylev:2003,bobylev:2008,bobylev:2002diff}, where in particular the uniqueness is obtained.

	    The physically relevant case of the hard sphere kernel $\lambda = 1$, $\gamma = 0$ was first considered in  \cite{Toscani:2000} in the unidimensional case. 
	    This work establishes the global existence of measure solutions with finite kinetic energy for this problem. The proof relies on a priori estimates of the solution to \eqref{eq:GranularGasesHomogeneous} in the Monge-Kantorovich-Wasserstein metrics $ W_2$.	    
	    This work also investigates the quasi-elastic $1-e^2 \sim \ve \to 0$ limit of the model, a nonlinear McNamara-Young-like friction equation.
	    This friction model was later investigated in \cite{Li:2004}, where its global wellposedness was shown in the space of measure of finite energy.
	    
		The tail behavior of the equilibrium solution to the granular gases equation with a thermal bath $\Delta_v f$ was investigated in many papers, the main ones being  \cite{ErnstBrito2002, bobylev:2004,gamba:2004}. They all proved the existence of non-gaussian, overpopulated tails for diffusively excited granular gases, namely:
			\begin{theorem}[From \cite{bobylev:2004} and \cite{ErnstBrito2002}]
			  \label{thm:TailBehavior}
			  Let $F(v) \geq 0$ for $v \in \RR^d$ be a solution to the stationary equation 
			  \[ \Q_\mathcal{I}(F,F) + \Delta_v F = 0\] 
			  with all polynomial moments in velocity. Then, 
			  \[ F(v) \sim_{|v| \to \infty} \exp \left (-|v|^{\alpha}\right ),\]
			  with $\alpha = 1$ in the Maxwellian molecules case and $\alpha = 3/2$ in the hard spheres case.
			\end{theorem} 			
		Indeed, the thermal bath gives an input of kinetic energy, preventing the appearance of trivial Dirac delta equilibria. The propagation of the Sobolev norms of the space dependent version of this equation was then established in \cite{gamba:2004}.
		It uses a careful estimate of the inelastic entropy production \eqref{inegEntropGranular}, and a fixed point argument for the existence and uniqueness of solutions.
			
		Finally, the work \cite{Mischler:20061} establishes the global well posedness of the granular gases equation without a thermal bath, for a general case of collision kernel (which contains \eqref{defHSkernelAnomalous}) and velocity dependent restitution coefficients:
			\begin{theorem}[Theorem 1.4 of \cite{Mischler:20061}]
				Let $0 \leq f_{in} \in L^1_3  \cap \in BV_4$. Then for any $T \in (0, T_c)$, where $T_c := \sup \, \{T > 0 : \mathcal{E}(f)(t) > 0, \, \forall \, t < T \}$ is the so-called \emph{blowup time}, there exists an unique nonnegative solution $f \in \mathcal C(0,T;L_2^1) \cap L^\infty (0,T; L_3^1)$ of \eqref{eq:GranularGasesHomogeneous}.
				It preserves mass and momentum, and converges in the weak-* topology of measures towards a Dirac delta.
			\end{theorem}
			Their proof relies on careful estimates of the collision operator $\Q_\mathcal{I}$ in Orlicz space (specially the $L \log L$ space of finite entropy measures).

        \begin{rem}
            The related (but still mostly open) problem of the propagation of chaos was considered in  \cite{MischlerMouhotWennberg:2011} for a very simplified inelastic collision operator with a thermal bath.
        \end{rem}

        \paragraph{\textbf{Cauchy problem in the space dependent setting.}}
        
		    The case of the space inhomogeneous setting\footnote{Physically more realistic, in part because of the spontaneous loss of space homogeneity that has been observed in \cite{Goldhirsch:1993}.} has been much less investigated. 
			
			The first result can be found in \cite{Benedetto:2002} for the model introduced in \cite{Toscani:2000} with a restitution coefficient given by \eqref{defRestiCoeffToscani}, in one dimension of space and velocity.
			This work establishes the existence and uniqueness of mild (perturbative) solutions, first for small  $L^1(dx \, dv)$ initial data, and then for compactly supported initial data.
			Their main argument is reminiscent from a work due to Bony in \cite{bony:1987} concerning discrete velocity approximation of the Boltzmann equation in dimension $1$.

			 The global existence of mild solutions in the general $\RR_x^3 \times \RR_v^3$ setting, for a large class of velocity-dependent restitution coefficient, but for initial data close to vacuum, was obtained in \cite{alonso:2009}.
			 The proof is based on a Kaniel-Shinbrot iteration on a very small functional space.
			 The stability in $L^1(\RR^3_x \times \RR^3_v)$ under the same assumptions was established  in \cite{Wu:2009}. 
			 Finally the existence and convergence to equilibrium in $\mathbb{T}_x^3 \times \RR_v^3$  for a weakly inhomogeneous granular gas\footnote{Namely, the initial condition is chosen with a lot of exponential moments in velocity, and close to a space homogeneous profile.} with a thermal bath was proved in \cite{Tristani:2013}, using a perturbative approach. 
		  
	\subsection{Large time behavior}
	
	    \smallskip
		\paragraph{\textbf{Macroscopic properties of the granular gases operator.}}
		  
		  Modeling-wise, the main microscopic difference between a granular gas and a perfect molecular gas is the dissipation of the kinetic energy. 
		  Using the weak form \eqref{defQGranWeak} among with the microscopic relations \eqref{eqMicroConsInel} of the inelastic collision operator, this yields
		  \begin{equation*}
		    \int_{\RR^d}  \Q_\mathcal{I}(f,f)(v) \begin{pmatrix} 1 \\ v \\ |v|^2 \end{pmatrix} dv \,=\,	 	  \begin{pmatrix} 0 \\ 0 \\ -D(f) \end{pmatrix},
		  \end{equation*}
		  where $D(f) \geq 0$ is the  \emph{energy dissipation} functional, which depends only on the collision kernel:
		  \begin{equation}
		    \label{def:DissipFunctional}
		    D(f) \,:=\, \int_{\mathbb{R}^d \times \mathbb{R}^d}f \, f_* \, \Delta \left(|v-v_*|, E(f) \right) dv \, dv_*.
		  \end{equation}
		  The quantity $\Delta  \left(|u|, E \right)$ is the so-called \emph{energy dissipation} rate, given using \eqref{eqMicroConsInel} by
		  \begin{equation}
		    \label{defDissipRate}
		    \Delta  \left(|u|, E \right) \,:=\, \frac{1-e^2}{4}\int_{\SSS^{d-1}} |u \cdot \omega|^2 \, B(|u|, \cos \theta, E ) \, d \omega \geq 0, \quad \forall e \in [0,1].
		  \end{equation}
		  
		  This dissipation of kinetic energy has a major consequence on the behavior of the solutions to the granular gases equation. Indeed, combined with the conservation of mass and momentum, it implies (at least formally) an \emph{explosive} behavior, namely convergence in the weak-* topology of solutions to \eqref{eq:GranularGases} towards Dirac deltas, centered in the mean momentum $\bm u$:
		  \[ f(t,\cdot) \rightharpoonup \delta_{v=\bm u}, \quad t \to \infty.\] 
		  
		  As for the entropy, it is not possible to obtain any entropy dissipation for this equation, in order to precise this large time behavior.
		  Indeed, as noticed in \cite{gamba:2004}, taking $\psi(v) = \log f(v)$ in \eqref{defQGranWeak} yields
		  \begin{align}
		    \int_{\RR^d} \Q_\mathcal{I}(f,f)(v) \,\log f(v)\,dv = & \
		      \frac{1}{2}\int_{\mathbb{R}^d \times \mathbb{R}^d \times \mathbb{S}^{d-1}} f_{*} \, f \log \left(\frac{f' \, f_*'}{f \,f_*} \right) B \, d\sigma \, dv \, dv_* \notag \\
		      = &  \ \frac{1}{2}\int_{\mathbb{R}^d \times \mathbb{R}^d \times \mathbb{S}^{d-1}} f_{*} \, f \left [\log \left(\frac{f' \, f_*'}{f \,f_*}\right ) - \frac{f' \, f_*'}{f \,f_*} + 1 \right] B \, d\sigma \, dv \, dv_* \notag \\ 
		      & + \frac{1}{2}\int_{\mathbb{R}^d \times \mathbb{R}^d \times \mathbb{S}^{d-1}} \left ( f_{*}' \, f' - f_{*} \, f \right ) B \, d\sigma \, dv \, dv_*. \label{inegEntropGranular}
		  \end{align}
		 The first term in \eqref{inegEntropGranular}, the elastic contribution, is nonpositive because $\log \lambda - \lambda + 1 \leq 0$ (this is Boltzmann's celebrated \emph{H Theorem}).
		 Nevertheless, the second term has no sign \emph{a priori}: it is 0 only in the elastic case (because of the involutive collisional transformation  $(v,v_*, \sigma) \to (v',v_*',\sigma)$). Boltzmann's entropy 
		 \[
		    \mathcal H(f) := \int_{\RR^{d_v}} f(v) \log f(v) \,dv
		 \]
		 is then not dissipated by the solution of the granular gases equation if $e <1$.
		 Some work has been done on that direction in the numerical side. Indeed, adding a drift term or a thermal bath in velocity can yield numerical entropy dissipation, as noticed in  \cite{SoriaMaynarMischlerMouhotReyTrizac:2015}.

        \paragraph{\textbf{Kinetic energy dissipation and the Haff's cooling Law.}}
        
          Let us assume in this subsection that the granular gas considered is space homogeneous, namely $f$ is solution to \eqref{eq:GranularGasesHomogeneous}.
		  Having no known entropy, one has then to use other macroscopic quantities to study the large time behavior of solutions to \eqref{eq:GranularGases}. Because of its explicit dissipation functional, kinetic energy is a good candidate for this. Moreover, being related to the variance, it allows to measure the concentration in velocity of the solution.
		  
		  In order to have an explicit bound for the energy dissipation, let us assume that the collision kernel is of the general type \eqref{defHSkernelAnomalous}.
		  Using polar coordinates, it is straightforward to compute the dissipation rate \eqref{defDissipRate}:
		  \begin{equation}
		    \label{eqDissipRate}
		    \Delta  \left(|u|, E\right) = b_1 \, \frac{1-e^2}{4} \,  |u|^{\lambda+2} E^\gamma,
		  \end{equation}		 
		  where thanks to \eqref{def:AngularCrossSection}
		  \[ b_1 = \left | \mathbb S^{d-2} \right | \int_{0}^{\pi} \cos^2(\theta) \sin^{d-3}(\theta) \, b(\cos(\theta)) \, d\theta < \infty.\]
		 
		  Using the conservation of mass and momentum, one can always assume that the initial condition is of unit mass and zero momentum. Plugging \eqref{eqDissipRate} into \eqref{def:DissipFunctional} then yields using Hölder and Jensen inequalities
		  \begin{align}
		    \frac{d}{dt} E(f)(t) & = - b_1 \, \frac{1-e^2}{4} \, E(f)^\gamma(t) \int_{\mathbb{R}^d \times \mathbb{R}^d} f \, f_* \, |v-v_*|^{\lambda+2} \,  dv \, dv_* \notag \\
		    & \leq - b_1 \, \frac{1-e^2}{4} \,  E(f)^\gamma(t) \int_{\RR^d} f(v) \, |v|^{\lambda+2} \, dv \notag \\
		    & \leq - b_1 \, \frac{1-e^2}{4} \,  E(f)^{1+\gamma + \lambda/2}(t). \label{inegEnerHS}
		  \end{align}
		  
		  
		  In particular, one will have the following large time behaviors: Setting $ C_e = b_1 \, \rho \, (1-e^2)/4$ and $\alpha := \gamma + 1/2 $, 
		  \begin{itemize}
		      \item Maxwellian pseudo-molecules ($\lambda = \gamma = 0$) decays exponentially fast towards the Dirac delta : 
		      \[  E(f)(t) = E\left (f_{in} \right ) e^{ - C_e \, t};\]
Notice that the inequality in \eqref{inegEnerHS} is an identity for this case.
		      
		      \item Hard spheres ($\lambda =1$, $\gamma =0$) exhibits the seminal quadratic \emph{Haff's cooling Law} \cite{haff:1983}:
		      \begin{equation} 
		        \label{eq:HaffLaw}
		        E(f)(t) \leq {\left (E\left (f_{in} \right )^{-1/2} + C_e \, t/2\right) ^{-2 }}.
		      \end{equation}
		      
		      \item Anomalous gases ($\gamma \neq 0$) exhibits more general behaviors: 
		  \begin{equation*}
		    E(f)(t) \leq \left\{ \begin{aligned}
		      & {\left (E\left (f_{in} \right )^\alpha + C_e \alpha \, t\right) ^{-\frac{1}{\alpha}} } && \text{if } \gamma >  {-1/2} \  (\alpha> 0, \ \text{ finite time extinction)}; \\
		      & E\left (f_{in} \right ) e^{ - C_e \, t}   && \text{if } \gamma = {-1/2}; \\
		      & \left ( E\left (f_{in} \right )^{\alpha} - C_e \alpha \, t\right)^{-\frac{1}{\alpha}}  && \text{if } \gamma < {-1/2} \ (\alpha < 0).
		     \end{aligned} \right.
		  \end{equation*}
		  \end{itemize}
		  All of these formal results have been proven to be rigorous and sharp, with explicit lower bounds, in \cite{bobylev:2000,BCT2} for the Maxwellian and hard sphere cases \cite{Mischler:20062}, and in \cite{rey:2011} for the anomalous cases. Extension to the viscoelastic case can be found in \cite{alonso:2009,alonso2014boltzmann}, where the energy is shown to behave as
		  \[
		    E(f)(t) \sim_{t \to \infty} C \, (1+t)^{-5/3}.
		  \]

		  All these papers share a common approach of proof, using the fact that the space homogeneous granular gases equation admits a self-similar behavior.
		  Hence, introducing some well chosen time-dependent scaling function $\omega$ and $\tau$, the distribution $f$ is written as
          \[ f(t,v) = \omega(t)^d g \left (\tau(t), \omega(t) \, v\right ), \]
		  to take into account the concentration in the velocity variables\footnote{One can see the velocity scaling function $\omega$ as the inverse of the variance of the distribution $f$. This scaling is then a continuous ``zoom'' on the blowup, and can be used to develop numerical methods for solving the full granular gases equation, see \cite{FR13}.}.
		  The rescaled function $g$ is then solution to the granular gases equation, with an \emph{anti-drift} term in velocity: 
		  \begin{equation*}
			  \partial_t g + \nabla_v \cdot (v \, g) = \mathcal Q_\mathcal{I}(g,g).
		  \end{equation*}
		  Using some regularity estimates of the gain term of $\Q_\mathcal{I}$ ``à la'' Lions/Bouchut-Desvillettes \cite{bouchut1998proof} and some new Povzner-like estimates \cite{AlCaGaMo:2012}, \cite{Mischler:20062}  then obtains a lower bound for the energy of $g$, yielding the generalized Haff's law by coming back to $f$.
		  
		  \begin{rem}
		    In the viscoelastic case, note that the rescaling in velocity induces a time dependency on the restitution coefficient, complicating the proof of the Haff's cooling law  \cite{alonso2013uniqueness}. It is also the case in the anomalous setting, where the rescaling function depends nonlinearly on the solution $f$.
		  \end{rem}

		  The question of the uniqueness, stability and exponential return to an universal equilibrium profile (\emph{hypocoercivity}, see \cite{Villani:2009Hypo}) of the self-similar solutions has then been fully addressed in the series of work \cite{mischler:20091, mischler:20092}, for a constant restitution coefficient, with and without a thermal bath. Extension of these results to the viscoelastic case has been done in \cite{alonso2014boltzmann}.
				
		\begin{rem}
		    These results are reminiscent of the works \cite{Benedetto:1998, CarrilloMcCann:2004,BCT1,BCT2,BC,Bolley:2012}, where the exponential convergence to equilibrium of the solution to the nonlocal granular media equation or the Maxwellian case has been shown in the $W_2$ or Fourier metrics.
	    \end{rem}

	\subsection{Compressible hydrodynamic limits}

        Let us consider in this subsection the following hyperbolic scaling of the granular gases equation:
        \begin{equation}
        \label{eq:GranularGasesScaled}
         \frac{\partial f_\ve}{\partial t} + v \cdot \nabla_{ x} f_\ve = \frac{1}{\ve} \Q_\mathcal{I}(f_\ve,f_\ve).
        \end{equation}	
        Determining the precise hyperbolic limit $\ve \to 0$ of equation \eqref{eq:GranularGasesScaled} is a fundamental, yet very difficult question.
        
        Indeed, for the elastic case $e=1$, one simply has to use the fact that the equilibria of the collision operator are at the thermodynamical equilibrium (gaussian distributions) and the conservation of mass, momentum and kinetic energy to obtain the classical compressible Euler-Fourier system \cite{CIP:94}.
        Because of the trivial Dirac equilibria, this question is more intricate for the true inelastic case.
            
        \smallskip
        \paragraph{\textbf{Pressureless Euler dynamics.}}
           	
            Adopting the same approach as in the elastic case, one can formally plug the ``equilibrium'' Dirac deltas in the pressure to obtain the following pressureless Euler system:
            \begin{equation} 
            	\label{sys:presurelessEuler}
            	\left\{ \begin{aligned}
            		& \frac{\partial \rho}{\partial t}  \,+\, \diverg  (\rho \, \bm{u}) \,=\, 0, 	\\
            		& \frac{\partial (\rho \bm{u})}{\partial t}  \,+\, \diverg  \left(\rho \, \bm{u} \otimes \bm{u} \right )\,=\, \bm{0}.
            	\end{aligned} \right.
            \end{equation}
            This system can describe various interesting physical situations, such as galactic clusters, but is notoriously difficult to study mathematically. 
            Its solution are in general ill-posed, as classical solutions cannot exists for large times and weak solutions are not unique. 
            
            In the unidimensional case, it is however possible to recover a well posed theory by imposing a semi-Lipschitz condition on $u$. This theory was introduced in \cite{BouchutJames:1999}, and later extended in \cite{boudin:2000} and \cite{HuangWang:2001}.
            We cite below the main result of \cite{HuangWang:2001}, where $M^1(\RR)$ denotes the space of Radon measures on $\RR$ and $L^2(\rho)$ for $\rho\geq 0$ in $M^1(\RR)$ denotes the space of functions which are square integrable against the density $\rho$.
            \begin{theorem} [From \cite{HuangWang:2001}]
            For any $\rho^0\geq 0$ in $M^1(\RR)$ and any $u^0\in L^2(\rho^0)$, there exists $\rho\in L^\infty(\RR_+, M^1(\RR))$ and $u\in L^\infty(\RR_+, L^2(\rho))$ solution to \eqref{sys:presurelessEuler} in the sense of distribution and satisfying the semi-Lipschitz Oleinik-type bound
            \begin{equation}
            u(t,x)-u(t,y)\leq \frac{x-y}{t},\quad \mbox{for}\ a.e.\;x> y.\label{Oleinik0}
            \end{equation}
            Moreover the solution is unique if $u^0$ is semi-Lipschitz or if the kinetic energy is continuous at $t=0$
            \[
            \int_\RR \rho(t,dx)\,|u(t,x)|^2\longrightarrow \int_\RR \rho^0(dx)\,|u^0(x)|^2,\qquad \mbox{as}\ t\rightarrow 0.
            \]\label{wellposedpressureless}
            \end{theorem}  
            The proof of Th. \ref{wellposedpressureless} is quite delicate, relying on duality solutions. For this reason, we only explain the rational behind the bound \eqref{Oleinik0}, which can be seen very simply from the discrete \emph{sticky particles} dynamics. We refer in particular to \cite{BrenierGrenier:1998} for the limit of this sticky particles dynamics as $N\rightarrow \infty$.
              
            Consider $N$ particles on the real line. We describe the $i^{\text{th}}$ particle at time $t> 0$ by its position $x_i(t)$ and its velocity $v_i(t)$. Since we are dealing with a one dimensional dynamics, we can always assume the particles to be initially ordered  
             \[
            x_1^{in} < x_2^{in} < \ldots < x_N^{in}. 
            \]
            The dynamics is characterized by the following properties     
                 \begin{enumerate}
             \item The particle $i$ moves with velocity $v_i(t)$: $\frac{d}{dt} x_i(t)=v_i(t)$.     
             \item The velocity of the $i^{\text{th}}$ particle is constant, as long as it does not collide with another particle: $v_i(t)$ is constant as long as $x_i(t)\neq x_j(t)$ for all $i\neq j$.
                   \item The velocity jumps when a collision occurs: if at time $t_0$ there exists $j \in \{1, \ldots, N \}$ such that $x_j(t_0 ) = x_i(t_0)$ and $x_j(t) \neq x_i(t)$ for any $t < t_0$, then all the particles with the same position take as new velocity the average of all the velocities
                     \begin{equation*}
                       v_i(t_0+) = \frac{1}{|{j | x_j(t_0 ) = x_i(t_0)}|}\sum_{j | x_j(t_0 ) = x_i(t_0)} v_j(t_0-).
                     \end{equation*}
                 \end{enumerate}
                Note in particular that particles having the same position at a given time will then move together at the same velocity. Hence, only a finite number of collisions can occur, as the particles aggregates.
            
            This property also leads to the Oleinik regularity. Consider any two particles $i$ and $j$ with $x_i(t)>x_j(t)$. Because they occupy different positions, they have never collided and hence $x_i(s)>x_j(s)$ for any $s\leq t$. If neither had undergone any collision then $x_i(0)=x_i(t)-v_i(t)\,t>x_j(0)=x_j(t)-v_j(t)\,t$ or
            \begin{equation}
            	     \label{eqMicroBound}
            	     \frac{\left (v_{i} - v_j\right )_+}{\left (x_{i}-x_j\right )_+} < \frac{1}{t}, 
            	   \end{equation}
                 where $x_+ := \max (x,0)$. It is straightforward to check that \eqref{eqMicroBound} still holds if particles $i$ and $j$ had some collisions
            between time $0$ and $t$.
            
            As one can see this bound is a purely dispersive estimate based on free transport and the exact equivalent of the traditional Oleinik regularization for Scalar Conservation Laws, see \cite{Oleinik}. It obviously leads to the semi-Lipschitz bound \eqref{Oleinik0} as $N\rightarrow \infty$.
	  
		    This result was extended to the one dimensional (in space and velocity) granular gases equation \eqref{eq:GranularGasesScaled} in \cite{JabinRey:2016}.
		    The basic idea of the proof of this work is to compare the granular gases dynamics to the pressureless gas system \eqref{sys:presurelessEuler}. The main difficulty is to show that $f_\ve$ becomes monokinetic at the limit (see also the very recent paper \cite{kang2019propagation}). This is intimately connected to the Oleinik property \eqref{Oleinik0}, just as this property is critical to pass to the limit from the discrete sticky particles dynamics.

            Unfortunately it is not possible to obtain \eqref{Oleinik0} directly. Contrary to the sticky particles dynamics, this bound cannot hold for any finite $\ve$ (or for any distribution that is not monokinetic). This is the reason why it is very delicate to obtain the pressureless gas system from kinetic equations (no matter how natural it may seem). Indeed we are only aware of one other such example in \cite{kang2015asymptotic}.

            One of the main contributions of \cite{JabinRey:2016} is a complete reworking of the Oleinik estimate, still based on dispersive properties of the free transport operator $v\,\partial_x $ but compatible with kinetic distributions that are not monokinetic, through the introduction of a new, global nonlinear energy.
            The main result in this paper is the following:
            \begin{theorem}[from \cite{JabinRey:2016}]
                Consider a sequence of weak solutions $f_\ve \in L^\infty([0,\ T],\ L^p(\RR^2))$ for some $p>2$ and with total mass $1$ to the granular gases Eq. \eqref{eq:GranularGasesScaled} such that all initial $v$-moments are uniformly bounded in $\ve$, some moment in $x$ is uniformly bounded, and $f^0_\ve$ is, uniformly in $\ve$, in some $L^p$ for $p>1$.
Then any weak-* limit $f$ of $f_\ve$ is monokinetic, $f=\rho(t,x)\,\delta(v-u(t,x))$ for $a.e.\ t$, where $\rho,\;u$ are a solution in the sense of distributions to the pressureless system \eqref{sys:presurelessEuler} while $u$ has the Oleinik property \eqref{Oleinik0}.
\end{theorem}		
            
        \paragraph{\textbf{Quasi-elastic limit.}}
        
        		    The physical community usually considers another approach, that is assuming that the granular gas is in a quasi-elastic $1-e^2\sim \varepsilon \rightarrow 0$ setting. This was first proposed in \cite{JR85,JR852}, using an approach very similar to the seminal Grad's 13 moments closure for rarefied gas dynamics. The difficulty of a hydrodynamic description of granular materials has been addressed in well reasoned terms in \cite{GoldhirschChaos99,LNP564}, and as already discussed in the introduction, the hydrodynamic equations obtained with the kinetic theory of granular gases have been shown to be insightful well beyond their supposed limit of validity, i.e., away from the quasi-elastic limit assumption with external sources of energy. In fact, assuming that solutions of the kinetic problem do not deviate from being Gaussians, one can then obtain in the hard sphere case the following quasi-elastic compressible Euler system
            \begin{equation} 
			\label{sys:compressibleEulerGranular}
				\left\{ \begin{aligned}
					& \frac{\partial \rho}{\partial t} \,+\, \diverg  (\rho \, \bm{u}) \,=\, 0, 
					\\
					& \frac{\partial (\rho \bm{u})} {\partial t} \,+\, \diverg  \left(\rho \,( \bm{u} \otimes \bm{u})  \,+\, \rho \, T \,{\rm\bf I}\right) \,=\, \bm{0}, 
					\\
					& \frac{\partial W}{ \partial t}     \,+\, \diverg  \left( \bm{u}  \,( W\,+\, \rho \, T ) \right) \,=\, -K \, \rho \, T^{3/2},
				\end{aligned} \right.
			\end{equation}
representing the evolution of number density of particles $\rho(t,x)$, velocity field~$\bm{u}(t,x)$ and the total energy $W$, which  is given by $W=\frac32 \rho T + \frac12 \rho |\bm{u}|^2=\rho \epsilon + \frac12 \rho |\bm{u}|^2$, with $T$ being the granular temperature (3D). Here $K = K(d, B )$ is an explicit nonnegative constant, see below. 
		    This is a compressible Euler-type system, which dissipates the kinetic energy thanks to its nonzero right hand side. The particular expression of this RHS allows, after integration in space, to recover the correct  Haff's cooling Law. 
		    
		    The assumption that the solutions are not far from Gaussians obviously degenerates in a free cooling granular gas at some point leading to the so-called clustering instability studied by means of \eqref{sys:compressibleEulerGranular}, see for instance \cite{CS} and the references therein. In fact, this assumption can be shown to be valid in the quasi-elastic limit, see \cite{mischler:20092} for a rigorous justification of this property. Physicists argue that this assumption is also generically true in practical experiments with external sources of energy such as the shock waves in granular  flows under gravity \cite{RBSS02,BMSS02}, clustering \cite{HM03,BSSP04,CS},  the Faraday instability for vibrating thin granular layers \cite{Faraday:1831,MeloUmbanhowarSwinney:1995,MeloUmbanhowarSwinney:1996,TsimringAranson:1997,UmbanhowarMeloSwinney:1998,BMSS02,Bougieetal04,carrillo:2008granular}, and many other applications, see \cite{carrillo:2008granular,johnson_gray:2011} and the references therein.
		    
		    Passing from the granular gases equation \eqref{eq:GranularGases} to \eqref{sys:compressibleEulerGranular} has not been established properly. It can still be done formally under the weak inelasticity hypothesis
		    $1-e^2 \sim \ve$, see \cite{Toscani:2004}. 
		    This particular scaling insures that the granular gases operator converges towards the elastic Boltzmann operator, as was shown rigorously in \cite{mischler:20091} in the space homogeneous setting. Moreover, it allows to characterize the equilibrium distribution of the limit operator, which is Gaussian.
		    
		    A first step towards a rigorous compressible hydrodynamic limit is available in \cite{rey:2013}, where the study of the spectrum of the heated granular gases operator $\Q_\I + \tau \Delta_v$, linearized with respect to the equilibrium described in \cite{mischler:20092}, is done. 
		    For small inelasticity $1-e^2 \sim 0$, this work provides a spectral decomposition, and more importantly the existence and computation of eigenvalue branches.
		    This work  generalizes the seminal paper \cite{ellis:1975} on the spectrum of the linearized Boltzmann operator in $L^2$ to the $L^1$ and inelastic setting.
		    
		    Other types of fluid limits (such as viscous limits) of the granular gases equation has been described in the review paper \cite{Dufty:2008} and in the recent survey \cite{Garzo2019} for many different physical regimes, but none has been rigorously established. To illustrate the kind of equations obtained through these procedure, we write the generalized Navier-Stokes compressible equations for granular media in conservative form, see \cite{carrillo:2008granular}, as
\begin{equation}\label{nscons}
\hspace*{-1cm}\begin{array}{l l}
&{\displaystyle \frac{\partial \rho}{ \partial t}} + \nabla\cdot (\rho \bm{u})  =  0   \vspace{6pt} \\ 
&{\displaystyle \frac{\partial (\rho \bm{u})}{ \partial t}} + \nabla\cdot \left[ \rho \, (\bm{u} \otimes \bm{u})\right]  =  \nabla\cdot \bm{P} + \rho \, \bm{F}  \vspace{6pt} \\ 
&{\displaystyle \frac{\partial W}{\partial t}} + \nabla\cdot\left[ \bm{u} W \right]  =  -\nabla\cdot q + \bm{P}:\bm{E}+ \bm{u} \cdot (\nabla\cdot \bm{P})- \gamma + \rho \bm{u} \cdot \bm{F}, 
\end{array}
\end{equation}
representing the evolution of number density of particles $\rho(t,x)$, velocity field~$\bm{u}(t,x)$ and the total energy $W$, which  is given by $W=\frac32 \rho T + \frac12 \rho |\bm{u}|^2=\rho \epsilon + \frac12 \rho |\bm{u}|^2$, with $T$ being the granular temperature (3D). The symbol $F$ stands for external forces applied to the system. The constitutive relations for the momentum and heat fluxes write, as usual,
\begin{equation*}
\hspace*{-8mm}P_{ij} = \left[ - p + \left( \lambda - \frac23 \mu \right) \sum_{i} E_{ii} \right] \delta_{ij} + 2\mu E_{ij}
\end{equation*}
for the stress tensor, with $E_{ij}=\frac12 \left( \frac{\partial U_i}{\partial x_j} + \frac{\partial U_j}{\partial x_i}\right)$. The thermal conductivity relates linearly the heat flux $q$ to the temperature gradient, $q = - \chi \nabla T$.

The equation of state is relevant here: we use the expression $G(\nu)= [1-(\nu/\nu_{max})^{\frac{4}{3}\nu_{max}}]^{-1} $ for the contact value of the pair correlation function $G(\nu)$, which accounts for high densities, and the equation of state $p = \rho T ( 1 + 2(1+e)\nu G(\nu)) $, where $\nu$ is the packing fraction and $e$ the constant restitution coefficient. Random close-packing is achieved in 3D at $\nu_{max}=0.65$;  we do not allow any packing fraction higher than 99.99\% of this value. The transport coefficients for constant restitution given in \cite{JR85,BSSS99} write,
\begin{equation*}
\hspace*{-8mm}\gamma =  \frac{12}{\sigma \sqrt{\pi}} (1-e^2) \rho T^{3/2} G(\nu),
\end{equation*}
for the cooling coefficient $\gamma$, which models energy dissipation through collisions. Other kinetic coefficients are the shear viscosity, 
\begin{equation*}
\hspace*{-8mm}\mu = \frac{\sqrt{\pi T}}{6} \rho \sigma 
\left[ \frac{5}{16G(\nu)} + 1 + \frac45 \left(1 + \frac{12}{\pi} \right) G(\nu) \right],
\end{equation*}
the  bulk viscosity, 
\begin{equation*}
\hspace*{-8mm}\lambda = \frac{8}{3\sqrt{\pi}} \rho \sigma \sqrt{T} G(\nu)
\end{equation*}
and the thermal conductivity 
\begin{equation*}
\hspace*{-8mm}\chi = \frac{15}{16}\sqrt{\pi T} \rho \sigma 
\left[ \frac{5}{24G(\nu)} + 1 + \frac{6G(\nu)}{5} \left(1 + \frac{32}{9\pi} \right) \right] ,
\end{equation*}
implemented directly with no further assumptions.  
		    
		    A first attempt to derive equations of compressible Navier-Stokes type was done in the paper \cite{carlen:2010} using singular perturbations  of the collision operator $\Q_\mathcal I$ and a central manifold approach inspired from \cite{fenichel:1979}. The fact is that transport coefficients for compressible Navier-Stokes like equations can be derived by moment closures under different assumptions and these equations are able to recover realistic phenomena in granular gases, see \cite{Garzo2019} for a very recent review. For instance, the Faraday instability for vertically oscillated granular layer was analysed in \cite{carrillo:2008granular,ACSGP} comparing molecular dynamics simulations, and different closures proposed in the literature. The conclusion, as the reader can check in the numerical experiments in \cite{online}, is that the qualitative behavior of the numerical experiments using equations \eqref{nscons} fully recovers the physical expected outcome. Thus, this is a physical validation of these kind of approximations that are totally opened for mathematicians to be derived in full rigor.

\section{Fourier spectral methods for the inelastic Boltzmann collision operator}

Due to the complexity of the inelastic Boltzmann collision operator $\Q_\mathcal{I}(f, f)$ (a five-fold nonlinear integral operator), numerical simulation of granular gases is challenging and mostly done at the particle level (direct simulation Monte Carlo method \cite{Bird} or its variants). Over the past decade, a class of deterministic numerical methods -- the Fourier-Galerkin spectral method -- has received a lot of popularity for its high accuracy and relatively low computational cost. The first attempt was made in \cite{NPT03} for the one-dimensional model. Later in \cite{FPT05, GT09, FR13}, both two and three dimensional cases were considered. Although the implementation details may differ, the essential ideas in these works are the same, that is, utilizing the translational invariance of the collision operator and convolution property of the Fourier basis to write the collision operator as a weighted convolution in the Fourier space. In this way, the $O(N^{3d})$ cost per evaluation of the collision operator in the Galerkin framework (since $\Q_\mathcal{I}(f, f)$ is quadratic) is readily reduced to $O(N^{2d})$, where $N$ is the number of basis used in each velocity dimension. Even though this reduction is dramatic compared to other spectral basis, numerical implementation of the ``direct'' Fourier spectral method is still computationally demanding; what makes it worse is that the method also requires $O(N^{2d})$ memory to store the precomputed weights, which quickly becomes a bottleneck as $N$ increases. Recently, a fast Fourier spectral method was introduced in \cite{HM19}, wherein the key idea is to shift some offline precomputed items to online computation so that the weighted convolution in the original method can be rendered into a few pure convolutions to be evaluated efficiently by the fast Fourier transform (FFT). As a result, both the computational complexity and the memory requirement in the direct Fourier method are reduced.

In this section, we briefly review the original direct Fourier spectral method proposed in \cite{FPT05} and then its fast version introduced in \cite{HM19}. To this end, let us introduce the following weak form of the inelastic Boltzmann collision operator $\Q_\mathcal{I}(f, f)$ using the $\sigma$-representation \eqref{eq:sigmaRep}: 
\begin{equation}
    \label{eq:weak_form}
  \int_{\mathbb{R}^d} \Q_\I(f,f)(v)\psi(v)\diff{v} = \int_{\mathbb{R}^d}\int_{\mathbb{R}^d}\int_{\mathbb{S}^{d-1}}B(|v-v_*|,\sigma\cdot\widehat{(v-v_*)})ff_*\left(\psi'-\psi\right)\diff{\sigma}\diff{v}\diff{v_*}\,,
\end{equation}
where $\sigma\cdot\widehat{(v-v_*)}=\cos\theta$ ($\hat{u}$ denotes the unit vector along $u$), and for simplicity we assume $B$ does not depend on the kinetic energy $E$ (compare with \eqref{defHSkernelAnomalous}). From the numerical point of view, this does not impose any limitation since the $E$ part can always be factored out from the integral sign.


\subsection{The direct Fourier spectral method}


We first perform a change of variable $v_* \rightarrow g=v-v_*$ in (\ref{eq:weak_form}) to obtain
\begin{equation*} 
  \int_{\mathbb{R}^d} \Q_\I(f,f)(v)\psi(v) \diff{v} = \int_{\mathbb{R}^d} \int_{\mathbb{R}^d} \int_{\mathbb{S}^{d-1}} B(|g|,\sigma\cdot \hat{g})f(v)f(v-g)  \left(\psi(v')-\psi(v)\right)\diff{\sigma} \diff{g} \diff{v},
\end{equation*}
where
\begin{equation*}
  v'=v-\frac{1+e}{4}(g-|g|\sigma).
\end{equation*}
We then assume that $f$ has a compact support: $\text{Supp}_v(f) \approx B_S$, where $B_S$ is a ball centered at the origin with radius $S$. Hence it suffices to truncate the infinite integral in $g$ to a larger ball $B_R$ with radius $R=2S$:
\begin{equation}\label{weak}
  \int_{\mathbb{R}^d} \Q_\I(f,f)(v)\psi(v)\diff{v} = \int_{\mathbb{R}^d} \int_{B_R} \int_{\SSS^{d-1}} B(|g|,\sigma\cdot \hat{g})f(v)f(v-g)  \left(\psi(v')-\psi(v)\right)\diff{\sigma} \diff{g} \diff{v}.
\end{equation}
Next we restrict $v$ to a cubic computational domain $D_L=[-L,L]^d$, and approximate $f$ by a truncated Fourier series:
\begin{equation*}
  f(v)\approx f_N(v)=\sum_{k=-\frac{N}{2}}^{\frac{N}{2}-1}\hat{f}_k e^{i\frac{\pi}{L}k\cdot v}, \quad \hat{f}_k=\frac{1}{(2L)^d}\int_{D_L}f(v)e^{-i\frac{\pi}{L}k\cdot v}\diff{v}.
\end{equation*}
Here $k=(k_1,\dots,k_d)$ is a multidimensional index, and $\sum_{k={-N/2}}^{N/2-1}:=\sum_{k_1,\dots,k_d=-N/2}^{N/2-1}$. The choice of $L$ should be chosen at least as $L\geq (3+\sqrt{2})S/2$ to avoid aliasing, see \cite{FPT05} for more details. Now substituting $f_N$ into (\ref{weak}) and choosing $\psi(v)=e^{-i\frac{\pi}{L}k\cdot v}$, we can obtain the $k$-th mode of the collision operator as
\begin{equation}\label{sum}
  \hat{Q}_k=\sum_{\substack{l,m=-\frac{N}{2}\\l+m=k}}^{\frac{N}{2}-1}G(l,m)\hat{f}_l\hat{f}_m,
\end{equation} 
where the weight $G(l,m)$ is given by
\begin{equation*}
  G(l,m)=\int_{B_R}e^{-i\frac{\pi}{L}m \cdot g}\left[\int_{\mathbb{S}^{d-1}}B(|g|,\sigma\cdot \hat{g})\left(e^{i\frac{\pi}{L}\frac{1+e}{4}(l+m)\cdot (g-|g|\sigma)}-1\right)\diff{\sigma}\right]\diff{g}.
\end{equation*}

In the original spectral method \cite{FPT05}, the weight $G(l,m)$ is precomputed and stored since it is independent of the solution $f$ which leads to a memory requirement of $O(N^{2d})$. During the online computation, the weighted sum (\ref{sum}) is directly evaluated whose complexity is $O(N^{2d})$.

\subsection{The fast Fourier spectral method}

To reduce the complexity of the direct spectral method as well as to alleviate its memory requirement, the key idea introduced in \cite{HM19} is to render the weighted convolution (\ref{sum}) into a pure convolution so that it can be computed efficiently by the FFT. One way to achieve this is through a low-rank approximation of $G(l,m)$, namely, 
\begin{equation}\label{lowrank}
  G(l,m)\approx\sum_{p=1}^{N_p}\alpha_p(l+m)\beta_p(m),
\end{equation}
where $\alpha_p$ and $\beta_p$ are some functions to be determined and the number of terms $N_p$ in the expansion is small. Then (\ref{sum}) becomes
\begin{equation} \label{approx}
  \hat{Q}_k\approx \sum_{p=1}^{N_p}\alpha_p(k)\sum_{\substack{l,m=-\frac{N}{2}\\l+m=k}}^{\frac{N}{2}-1}\hat{f}_l  \left(\beta_p(m)\hat{f}_m\right),
\end{equation} 
where the inner summation is a pure convolution of two functions $\hat{f}_l$ and $\beta_p(m)\hat{f}_m$. Hence the total complexity to evaluate $\hat{Q}_k$ (for all $k$) is brought down to $O(N_pN^d\log N)$, i.e., a few number of FFTs.

Specifically, we first split $G(l, m)$ into a gain term and a loss term:
\begin{equation*}
  G(l,m) = G_\text{gain}(l,m) - G_\text{loss}(m),
\end{equation*} 
where
\begin{equation*}
  G_\text{gain}(l,m):= \int_{B_R}e^{-i\frac{\pi}{L}m \cdot g}\left[\int_{\mathbb{S}^{d-1}}B(|g|,\sigma\cdot \hat{g})e^{i\frac{\pi}{L}\frac{1+e}{4}(l+m)\cdot (g-|g|\sigma)}\diff{\sigma}\right]\diff{g},
\end{equation*}
and
\begin{equation*}
  G_\text{loss}(m):=\int_{B_R}e^{-i\frac{\pi}{L}m \cdot g}\left[\int_{\SSS^{d-1}}B(|g|,\sigma\cdot \hat{g})\diff{\sigma}\right]\diff{g}.
\end{equation*}
Note that the loss term is readily a function of $m$, hence no approximation/decomposition is actually needed. This suggests to evaluate the loss term of the collision operator by a precomputation of $G_\text{loss}(m)$ and then compute
\begin{equation*}
 \hat{Q}^-_k=\sum_{\substack{l,m=-\frac{N}{2}\\l+m=k}}^{\frac{N}{2}-1}\hat{f}_l\left(G(m)\hat{f}_m\right)
\end{equation*}
by FFT. For the gain term, to get a decomposition of form (\ref{lowrank}), we introduce a quadrature rule to discretize $g$, then $G_\text{gain}(l, m)$ can be approximated as
\begin{equation} \label{lowrank1}
  G_\text{gain}(l,m)\approx \sum_{\rho, \hat{g}}w_{\rho}w_{\hat{g}}\, \rho^{d-1} e^{-i\frac{\pi}{L}\rho\, m \cdot \hat{g}} F(l+m,\rho,\hat{g}),
\end{equation}
where $\rho:=|g|\in [0,R]$ is the radial part of $g$ and $\hat{g}\in \SSS^{d-1}$ is the angular part, and $w_{\rho}$ and $w_{\hat{g}}$ are the corresponding quadrature weights. The function $F$ is given by
\begin{equation}\label{FF}
  F(l+m,\rho,\hat{g}):=\int_{\SSS^{d-1}}B(\rho,\sigma\cdot \hat{g})e^{i\frac{\pi}{L}\rho \frac{1+e}{4}(l+m)\cdot (\hat{g}-\sigma)}\diff{\sigma}.
\end{equation}
With this approximation, the gain term of the collision operator can be evaluated as
\begin{equation*}
\hat{Q}_k^+\approx\sum_{\rho,\hat{g}}w_{\rho}w_{\hat{g}}\,\rho^{d-1}
      \,F(k,\rho,\hat{g})\sum_{\substack{l,m=-\frac{N}{2}\\l+m=k}}^{
        \frac{N}{2}-1}\hat{f}_l \left[e^{-i\frac{\pi}{L}\rho\, m \cdot
        \hat{g}}\hat{f}_m\right],
\end{equation*}
which is in the same form as explained in (\ref{approx}).

As for the quadratures, the radial direction $\rho$ can be approximated by the Gauss-Legendre quadrature. Since the integrand in (\ref{lowrank1}) is oscillatory on the scale of $O(N)$, the number of quadrature points needed for $\rho$ should be $O(N)$. The angular direction in $2$D can be discretized using simple rectangular rule which is expected to yield spectral accuracy due to the periodicity. While in 3D, we choose to use the spherical design \cite{Womersley2018} which is the near optimal quadrature on the sphere. 

To summarize, the total complexity to evaluate $\hat{Q}_k$ is $O(MN^{d+1}\log N)$, where $M$ is the number of points used on $\SSS^{d-1}$ and $M\ll N^{d-1}$. Furthermore, the only quantity that needs to be precomputed and stored is (\ref{FF}), which in the worst scenario only requires $O(MN^{d+1})$ memory.


\section{Numerical experiments and results}

The accuracy and efficiency of the fast spectral method has been validated in \cite{HM19}. In this section, we perform some additional tests to demonstrate the potential of the method in predicting some mathematical theories. We also introduce a GPU implementation of the method that significantly improves the CPU version in \cite{HM19}. This is critical especially for long time simulation.

We consider the following spatially homogeneous equation with a heat bath:
\begin{equation}\label{eqn:homo_heat}
  \partial_t f = \Q_\I(f, f) + \tau \Delta_v f,
\end{equation}
where $\tau$ is the parameter describing the strength of the heat bath. Notice that 
it is not necessarily related to the inelasticity parameter $e$, contrarily to \textit{e.g.} \cite{mischler:20092}.
The heat bath $\Delta_v f$ will also be discretized using the Fourier spectral method and Runge-Kutta method is used for time marching.

For the collision operator, we consider the simplified variable hard sphere kernel
\begin{equation}\label{VHS}
  B(|g|,\sigma\cdot\hat{g},E)=C_{\lambda}|g|^{\lambda}, \quad 0\leq \lambda\leq 1,
\end{equation}
where $C_{\lambda}>0$ is some constant (namely \eqref{defHSkernelAnomalous} with $b = C_\lambda$ and $\gamma = 0$.

For Maxwell molecules, given the initial condition $f_0(v)$ whose macroscopic
quantities are $\rho_0$, $u_0$ and $T_0$, the density and velocity are conserved
so $\rho(t)=\rho_0$, $u(t) = u_0$ and the temperature will evolve as
\begin{equation}\label{eq:ext_temp}
  T(t) = \left(T_0-\frac{8\tau}{1-e^2}\right)\exp{\left(-\frac{\rho_0(1-e^2)}{4}
  t\right)}+\frac{8\tau}{1-e^2},
\end{equation}
We could see
\begin{equation*}
  \lim\limits_{t\rightarrow\infty}T(t) = \frac{8\tau}{1-e^2}.
\end{equation*}
As in~\cite{HM19}, this analytical formula of temperature works as the
reference solution to ensure the correctness of the numerical solution. 

\subsection{GPU parallelized implementation}

From the implementation perspective, we dramatically improve the efficiency of
the fast spectral method by using GPU via Nivida's CUDA. GPU-parallelized
implementations are run on $2$
Intel\textsuperscript{\textregistered} Xeon\textsuperscript{\textregistered}
Silver 4110 2.10 GHz CPUs with $4$ NVIDIA Geforce GTX 2080 Ti (Turing) GPUs
accompanying CUDA driver 10.0 and CUDA runtime 10.0. The operating system used
is 64-bit Unbuntu 18.04. The CPU has 8 cores, 16 threads with max turbo
frequency 3.00 Ghz equiped with 128GB DDR4 REG ECC memory. The GPU has 4352 CUDA
cores, 11GB device memory. The algorithm has been written in python with
packages Numpy and Scipy. The CPU implementation is based on PyFFTW which is a
python wrapper of the C library FFTW. The GPU implementation is based on CuPy
which is an implementation of NumPy-compatible multi-dimensional array on CUDA.
As shown in Table~\ref{GPU_boost}, GPU version is up to $15$ times faster than
CPU version depending on different $N$s.
\begin{table}
	\centering
  \begin{tabu} to 0.6\linewidth {X[1, c] X[3, c] X[3, c]}
    \toprule
    $N$ & CPU & GPU \\
    \midrule
    8 & 7.68ms & 5.89ms \\
    16 & 61.2ms & 5.97ms \\
    32 & 546ms & 12.1ms \\
    64 & 5.38s & 109ms \\
    \bottomrule
  \end{tabu}
  \caption{Average running time per evaluation of the collision operator in
  $3$D. Comparison between the CPU and GPU-parallelized implementation for
  various $N$s ($\#$ of Fourier basis in each velocity dimension) and fixed $N_{\rho}=30$, $M_{\text{sph}}=32$ ($\#$ of quadrature points used in radial and spherical direction, respectively) with $2$
  Xeon\textsuperscript{\textregistered} Silver 4110 2.10 GHz CPUs with $4$
  NVIDIA Geforce GTX 2080 Ti (Turing) GPUs.}
	\label{GPU_boost}
\end{table}

\subsection{Numerical results}
\smallskip
\paragraph{\textbf{Test 1. Validation of exponential convergence to the equilibrium in 2D}}

As a first test, our goal is to verify the theoretical result
of the exponential convergence to the equilibrium. We consider the the Maxwell molecule by taking the
collision kernel as 
\begin{equation}\label{eq:maxw_kern}
  B(|g|, \sigma\cdot \hat{g},E) = C_0 = \frac{1}{2\pi}.
\end{equation}
The initial condition is chosen as
\begin{equation}\label{eq:maxw_init}
  f_0(v) = \frac{\rho_0}{2\pi T_0}e^{-\frac{(v - u_0)^2}{2T_0}},
\end{equation}
with $\rho_0=1, u_0=(0,0)$ and $T_0 = 8$. 

\begin{figure}
  \centering
  \includegraphics[width=\textwidth]{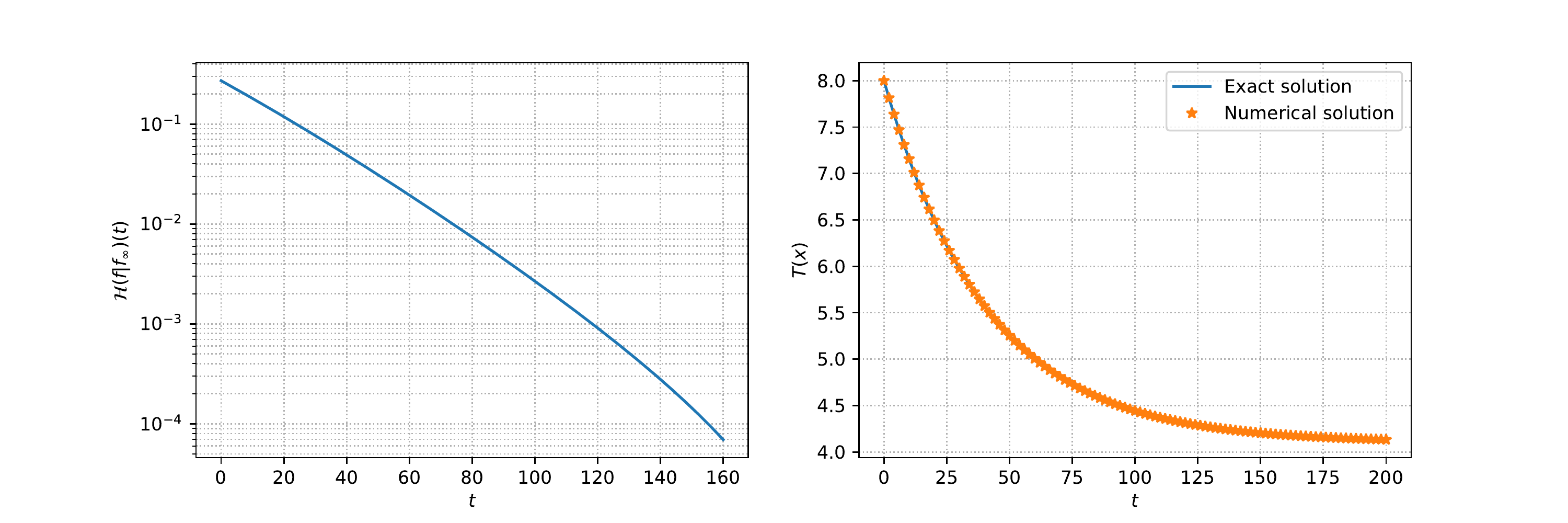}
  \caption{\textbf{Test 1.} Convergence to equilibrium. $e=0.95$ with heat bath $\tau=0.05$.
    Initial data is the Maxwellian~\eqref{eq:maxw_init}. Left: Semi-log plot of
    the relative entropy of $f$ and $f_\infty=f(t=100,v)$. Right: numerical temperature
    (orange dots) with exact temperature (blue line)
    as~\eqref{eq:ext_temp}. Numerical parameters: $N_v^2 = 64\times 64$, $N_{\rho}=32$,
    $M_{\hat{g}} = 16$, $R = 20$, $L = 5(3+\sqrt{2})$ and $\Delta t=0.01$.}
  \label{fig:convergence_095}
\end{figure}

\begin{figure}
  \centering
  \includegraphics[width=\textwidth]{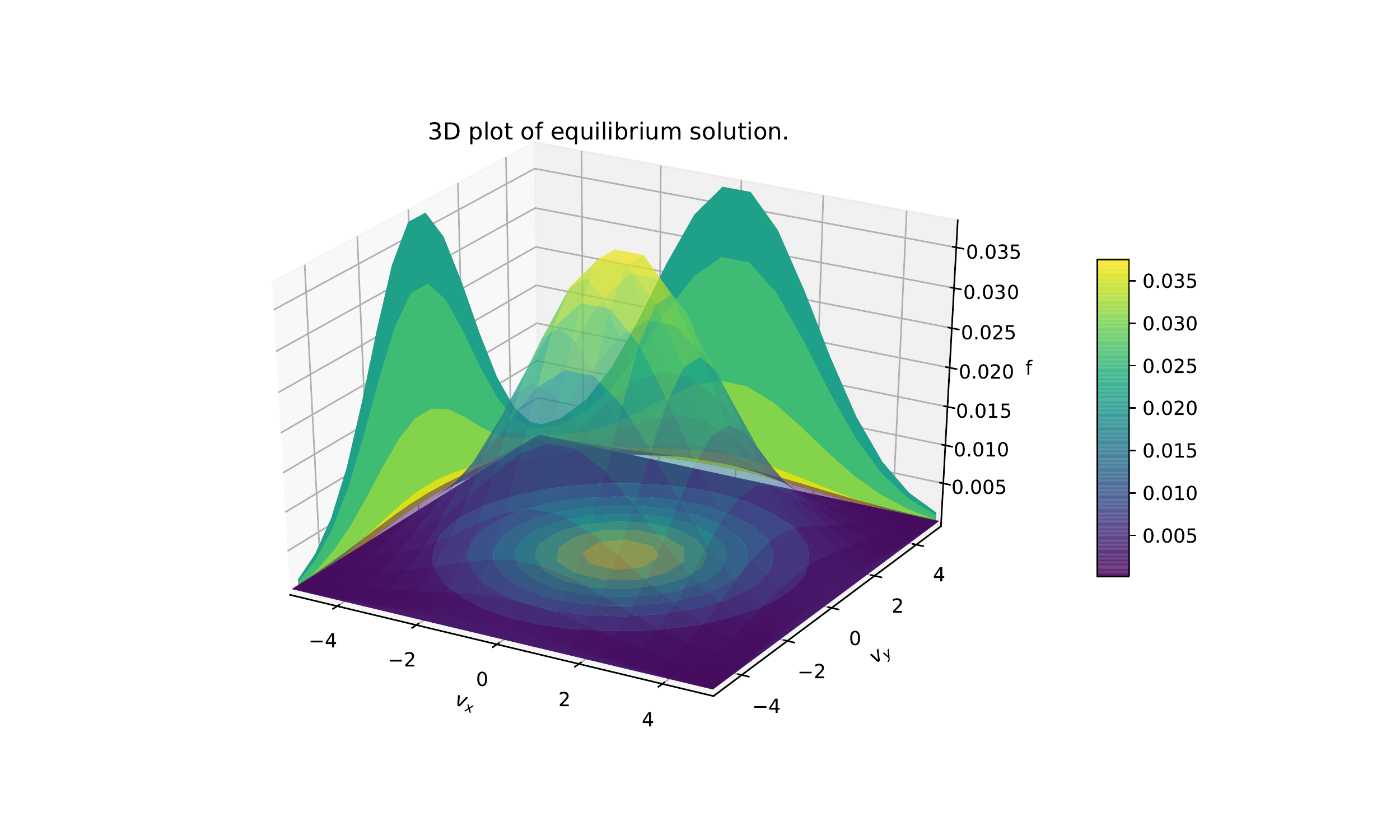}
  \caption{\textbf{Test 1.} The equilibrium profile of $e=0.95$ with heat bath $\tau=0.05$,
  initial data is the Maxwellian~\eqref{eq:maxw_init}. Numerical parameters:
  $N_v^2 = 64\times 64$, $N_{\rho}=32$, $M_{\hat{g}} = 16$, $R = 20$,
  $L = 5(3+\sqrt{2})$ and $\Delta t=0.01$.}
  \label{fig:equi_prof_095}
\end{figure}

The theoretical results from \cite{bobylev:2004} states that
if $e = 1-\tau$, then there exists a unique equilibrium solution $f_\infty$
to~\eqref{eqn:homo_heat} and one has
\begin{equation}\label{eq:exp_conv_rH}
  \mathcal{H}(f|f_\infty)(t)\sim O(e^{-\lambda t}),
\end{equation}
where the relative entropy is defined as
\begin{equation*}
  \mathcal{H}(f|g) = \int f\ln\left(\frac{f}{g}\right)\diff{v}\,.
\end{equation*}
In order to confirm this result, we first choose the following physical
parameters
\begin{equation*}
  e = 0.95, \quad \tau = 1-e = 0.05\,.
\end{equation*}
The numerical parameters we use to compute the $2$D collision operator are
\begin{equation*}
  N_v^2 = 64\times 64, \quad N_{\rho}=32, \quad M_{\hat{g}} = 16,\quad R = 20,
  \quad L = 5(3+\sqrt{2}).
\end{equation*}
For time integrator, a $4$-th order Runge-Kutta method is used with
$\Delta t = 0.01$ and we compute sufficient long time to get $f_\infty$ and in both cases the $\ell^2$ difference of solutions between the last $2$ time steps are of $O(10^{-12})$. The results are shown in Fig.~\ref{fig:convergence_095}
where one can observe the exponential convergence of relative
entropy~\eqref{eq:exp_conv_rH} from the left figure. In the right figure we find
a perfect match of the temperature evolution between numerical solution and the
analytical solution~\eqref{eq:ext_temp}.
We also plot the profile of the equilibrium solution $f_\infty$ in
Fig.~\ref{fig:equi_prof_095} which shows a Gaussian-like density
function.

\begin{figure}
  \centering
  \includegraphics[width=\textwidth]{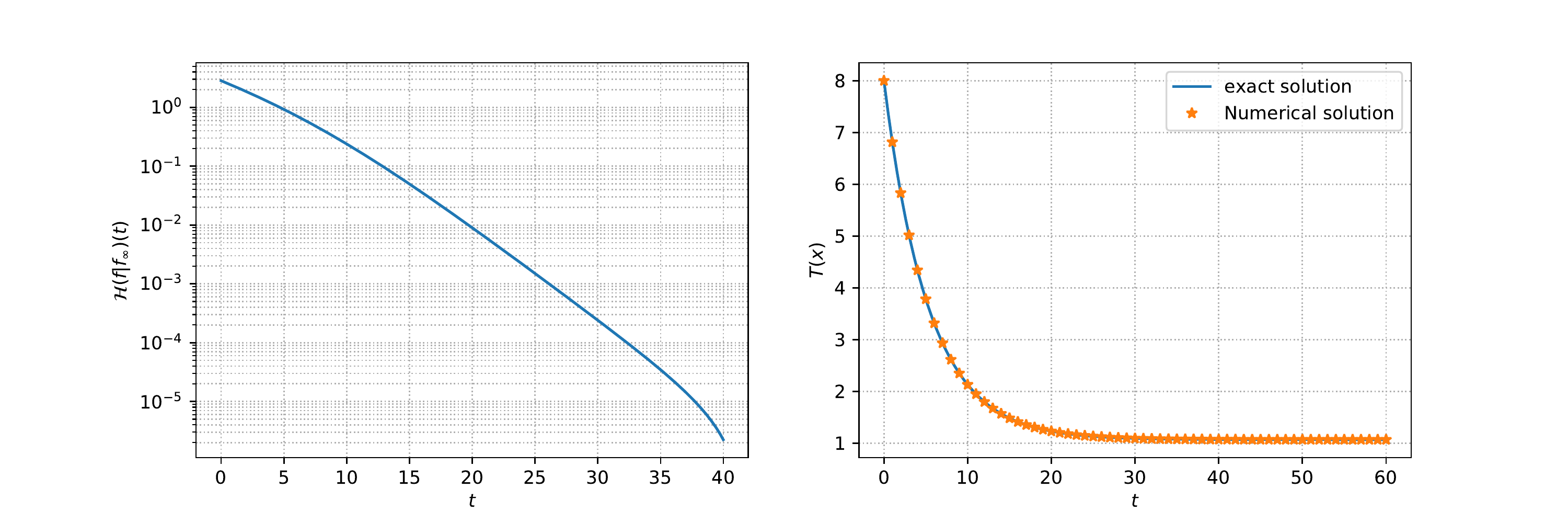}
  \caption{\textbf{Test 1.} Convergence to equilibrium. $e=0.5$ with heat bath $\tau=0.1$,
    initial data is the Maxwellian~\eqref{eq:maxw_init}. Left: Semi-log plot of
    the relative entropy of $f$ and $f_\infty=f(t=55,v)$. Right: numerical temperature
    (orange dots) with exact temperature (blue line)
    as~\eqref{eq:ext_temp}. Numerical parameters: $N_v^2 = 64\times 64$, $N_{\rho}=32$,
    $M_{\hat{g}} = 16$, $R = 20$, $L = 5(3+\sqrt{2})$ and $\Delta t=0.01$.}
  \label{fig:convergence_05_maxw}
\end{figure}

\begin{figure}
  \centering
  \includegraphics[width=\textwidth]{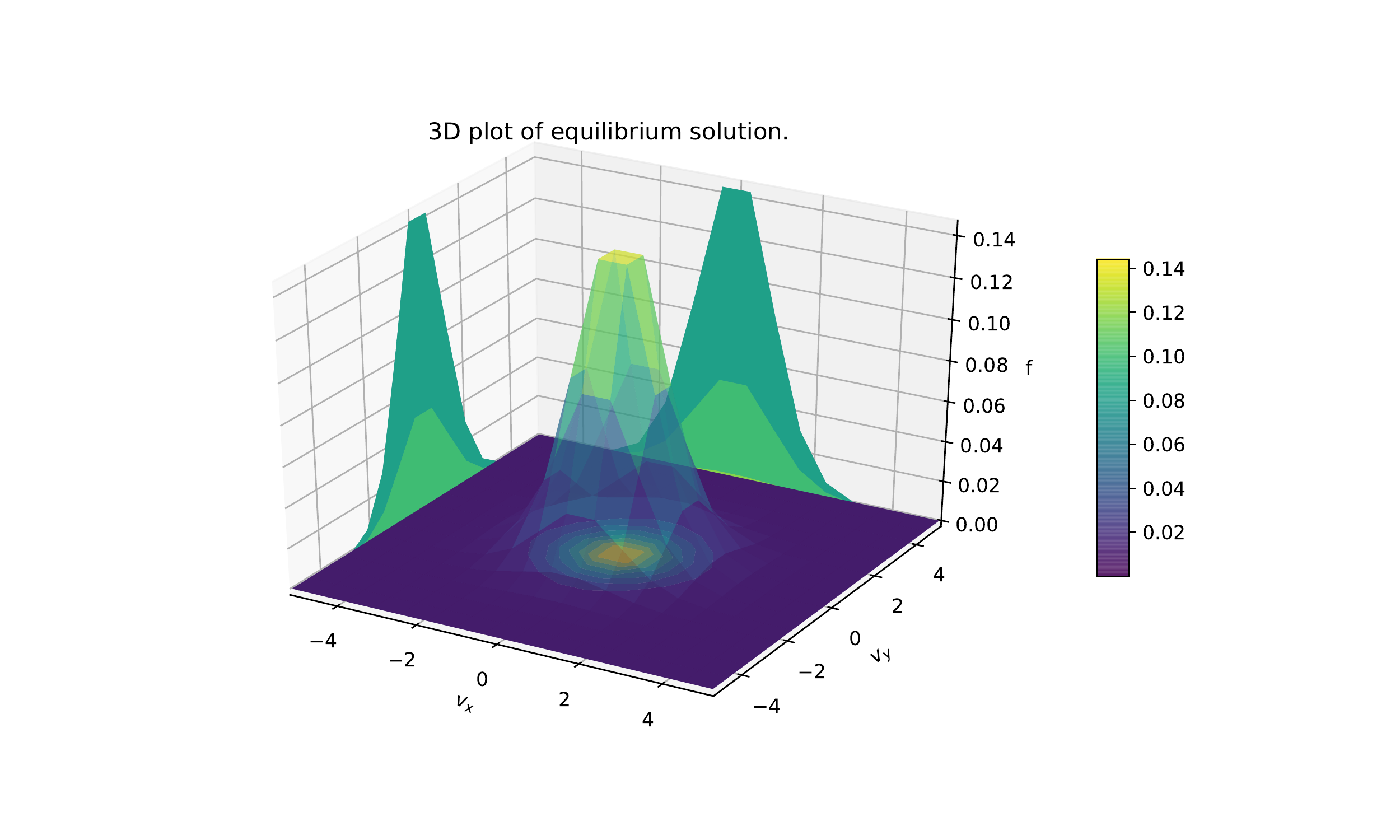}
  \caption{\textbf{Test 1.} The equilibrium profile of $e=0.5$ with heat bath $\tau=0.1$,
    initial data is the Maxwellian~\eqref{eq:maxw_init}. Numerical parameters:
    $N_v^2 = 64\times 64$, $N_{\rho}=32$, $M_{\hat{g}} = 16$, $R = 20$,
    $L = 5(3+\sqrt{2})$ and $\Delta t=0.01$.}
    \label{fig:equi_prof_05_maxw}
\end{figure}

Although the exponential convergence result is only available for the case that $e = 1 - \tau$, we expect something similar happens even when $e$ and $\tau$ are not related. In order to investigate this, we choose
\begin{equation*}
    e = 0.5, \quad\tau = 0.1,
\end{equation*}
and perform the test using the same initial data~\eqref{eq:maxw_init}. The
results are shown in Fig.~\ref{fig:convergence_05_maxw} and
Fig.~\ref{fig:equi_prof_05_maxw} where we indeed observe the same 
exponential convergence to equilibrium.

\begin{figure}
  \centering
  \includegraphics[width=\textwidth]{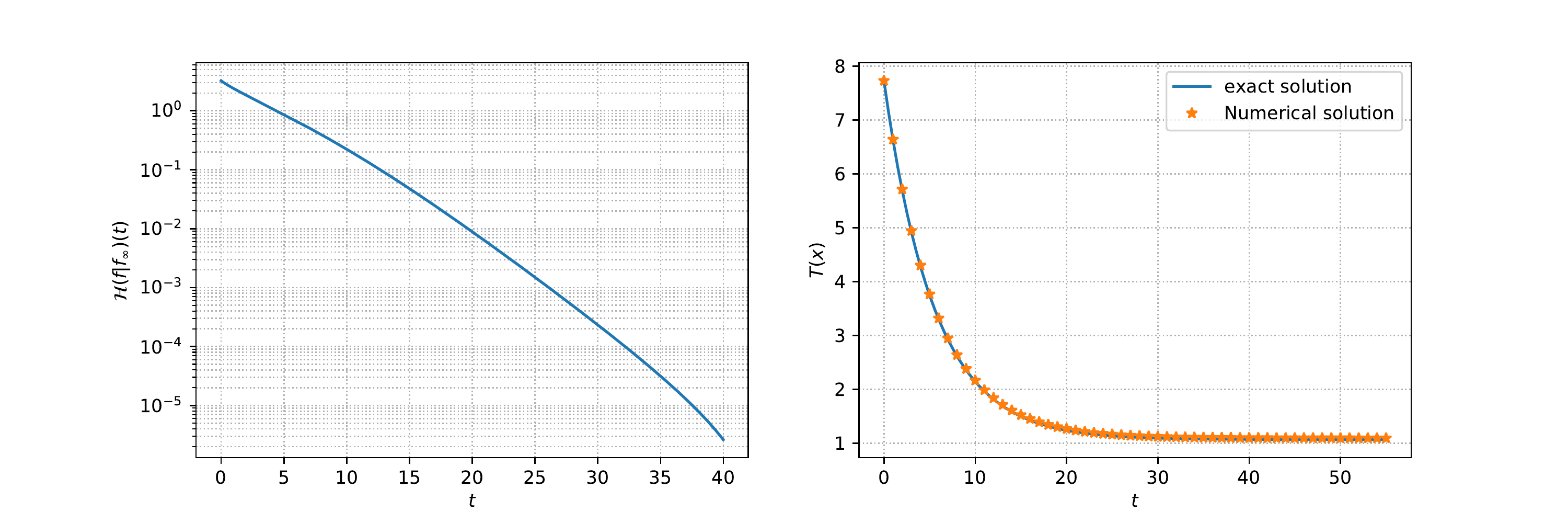}
  \caption{\textbf{Test 1.} Convergence to equilibrium. $e=0.5$ with heat bath $\tau=0.1$,
    initial data is the flat function~\eqref{eq:flat_init}. Left: Semi-log plot
    of the relative entropy of $f$ and $f_\infty=f(55,v)$. Right: numerical temperature
    (orange dots) with exact temperature (blue line)
    as~\eqref{eq:ext_temp}. Numerical parameters: $N_v^2 = 64\times 64$, $N_{\rho}=32$,
    $N_{\hat{g}} = 16$, $R = 20$, $L = 5(3+\sqrt{2})$ and $\Delta t=0.01$.}
  \label{fig:convergence_05_flat}
\end{figure}

\begin{figure}
  \centering
  \includegraphics[width=\textwidth]{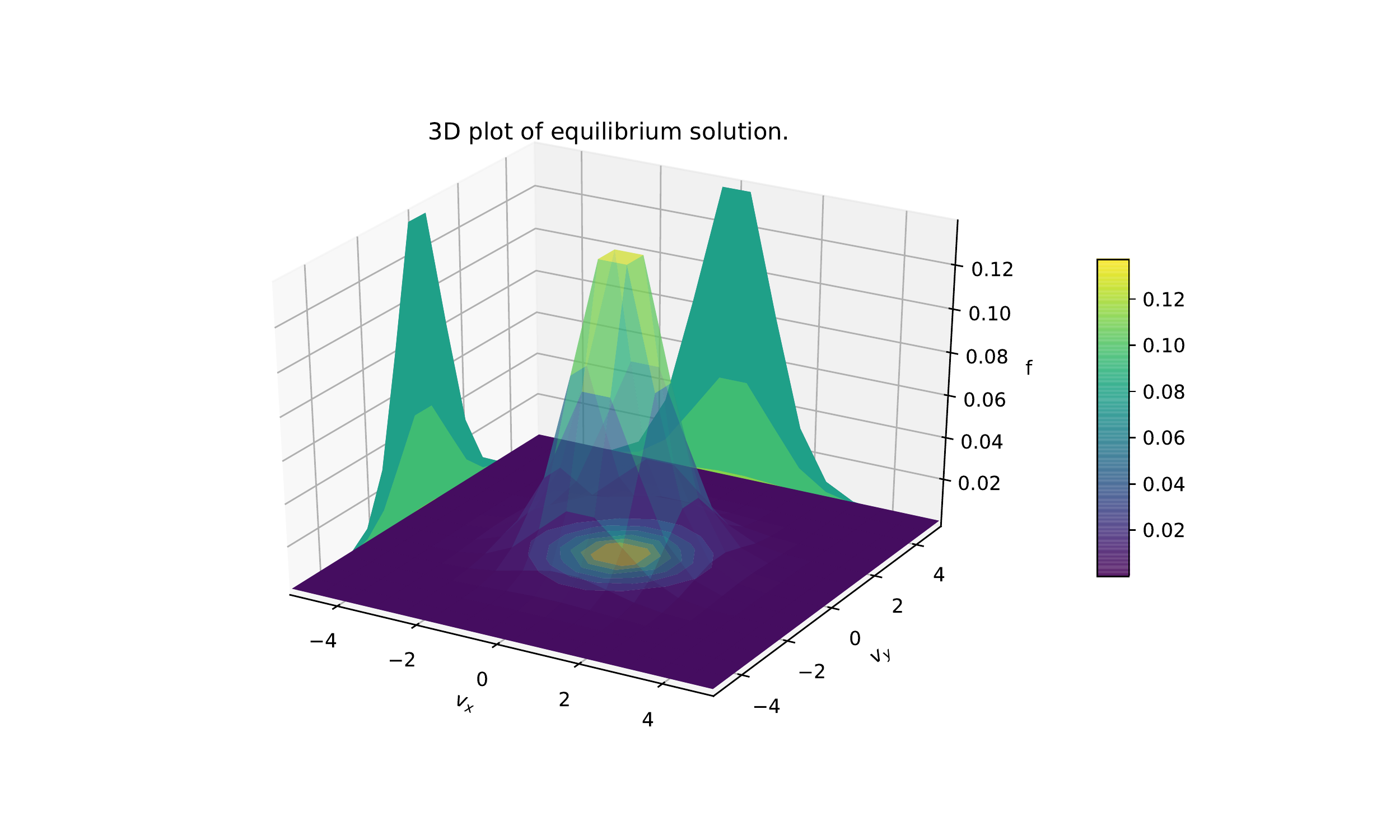}
  \caption{\textbf{Test 1.} The equilibrium profile of $e=0.5$ with heat bath $\tau=0.1$,
    initial data is the flat function~\eqref{eq:flat_init}. Numerical
    parameters: $N_v^2 = 64\times 64$, $N_{\rho}=32$, $M_{\hat{g}} = 16$, $R = 20$,
    $L = 5(3+\sqrt{2})$ and $\Delta t=0.01$.}
    \label{fig:equi_prof_05_flat}
\end{figure}

To confirm that we could always get Gaussian-like equilibrium
solution regardless of the initial data, we also consider the following initial condition
\begin{equation}\label{eq:flat_init}
  f_0(v) =
  \begin{cases}
    \dfrac{1}{4w_0^2},&\quad\text{for}\quad v\in [-w_0, w_0]\times[-w_0, w_0],\\
    0,                &\quad\text{otherwise},
  \end{cases}
\end{equation}
with $w_0 = 2\sqrt{6}$ such that $\rho_0=1, u_0=(0,0)$ and $T_0 = 8$. With
restitution coefficient $e=0.5$ and heat bath $\tau=0.1$, the results are shown in
Fig.~\ref{fig:convergence_05_flat} and Fig.~\ref{fig:equi_prof_05_flat}.

\smallskip
\paragraph{\textbf{Test 2. Investigation of tail behavior of the equilibrium in 2D}}

\begin{figure}
  \centering
  \includegraphics[width=\textwidth]{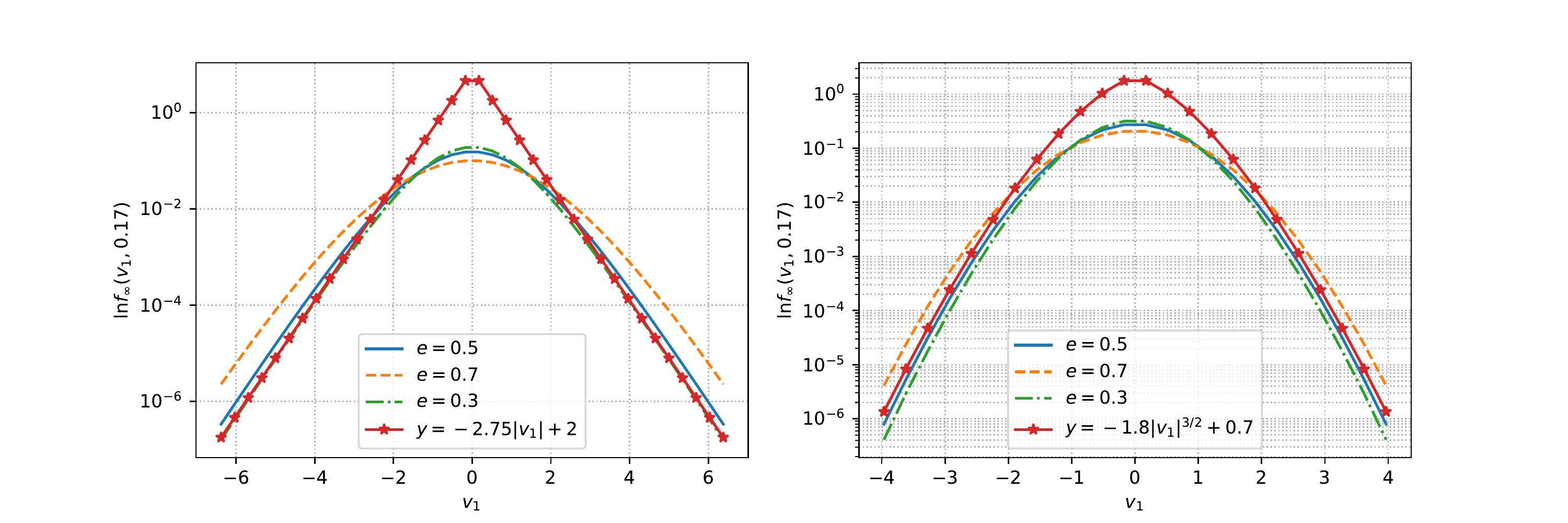}
  \caption{\textbf{Test 2.} The equilibrium profile of $e=0.3, 0.5, 0.7$ with heat bath
  $\tau=0.1$, initial data is the flat function~\eqref{eq:flat_init}. Left:
  Semi-log plot of $f_\infty(v_1, 0.17)=f(t=55, v_1, 0.17)$ for Maxwell molecules. Right: Semi-log
  plot of $f_\infty(v_1, 0.17) = f(t=55, v_1, 0.17)$ for hard spheres. The red lines are the reference
  profiles. Numerical parameters: $N_v^2 = 128\times 128$, $N_{\rho}=32$,
  $M_{\hat{g}} = 16$, $R = 20$, $L = 5(3+\sqrt{2})$ and $\Delta t=0.01$.}
  \label{fig:High_energy_tail}
\end{figure}

We now compare the different tail behaviors of the equilibrium solutions for the  Maxwell molecules collision kernel \eqref{eq:maxw_kern} and for the  hard spheres collision kernel \[
B(|g|,\sigma\cdot\hat{g})=|g|/(2\pi)
\]
in $2$D. To see the tail we need
higher resolution in velocity space so the velocity mesh is increased to
$N_v^2=128\times 128$. We plot the profile in $v_x$ ($v_1$) by choosing a fixed
$v_y$ ($v_2$) for different $e$s ($0.3$, $0.5$ and $0.7$). From
Fig.~\ref{fig:High_energy_tail}, we see that the numerical scheme generates overpopulated equilibrium tails: the  Maxwell molecules case behaves like
\begin{equation*}
  f(v, t=\infty)\sim e^{-\alpha |v|}\,,
\end{equation*}
and the hard spheresones behaves like
\begin{equation*}
  f(v, t=\infty)\sim e^{-\alpha |v|^{3/2}}\,.
\end{equation*}
These results corresponds accurately to what was predicted theoretically in \cite{ErnstBrito2002,bobylev:2004} (summarized in Theorem \ref{thm:TailBehavior}).

\smallskip
\paragraph{\textbf{Test 3. 3D hard sphere.}}

\begin{figure}
  \centering
  \includegraphics[width = .48\textwidth]{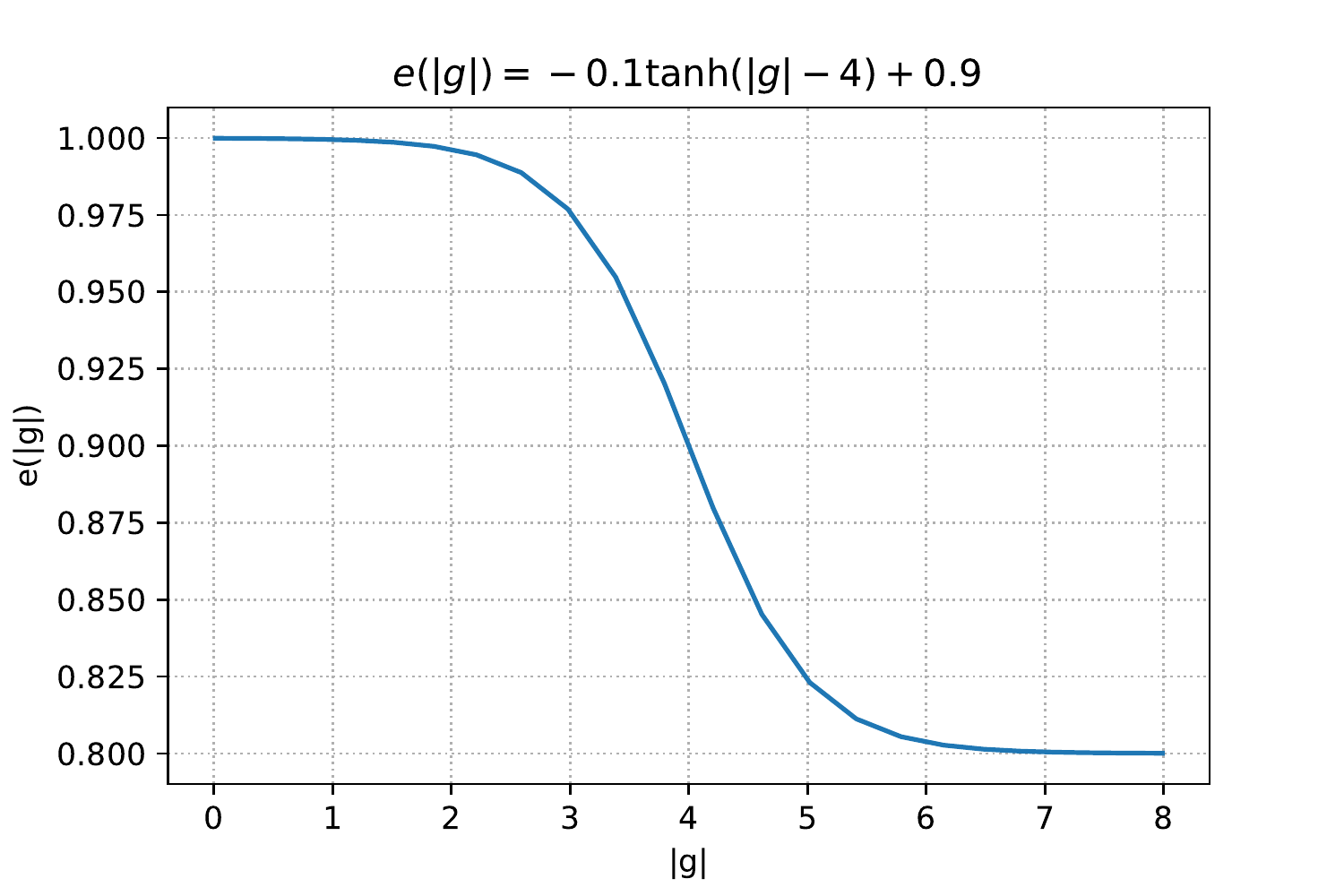}
  \includegraphics[width = .48\textwidth]{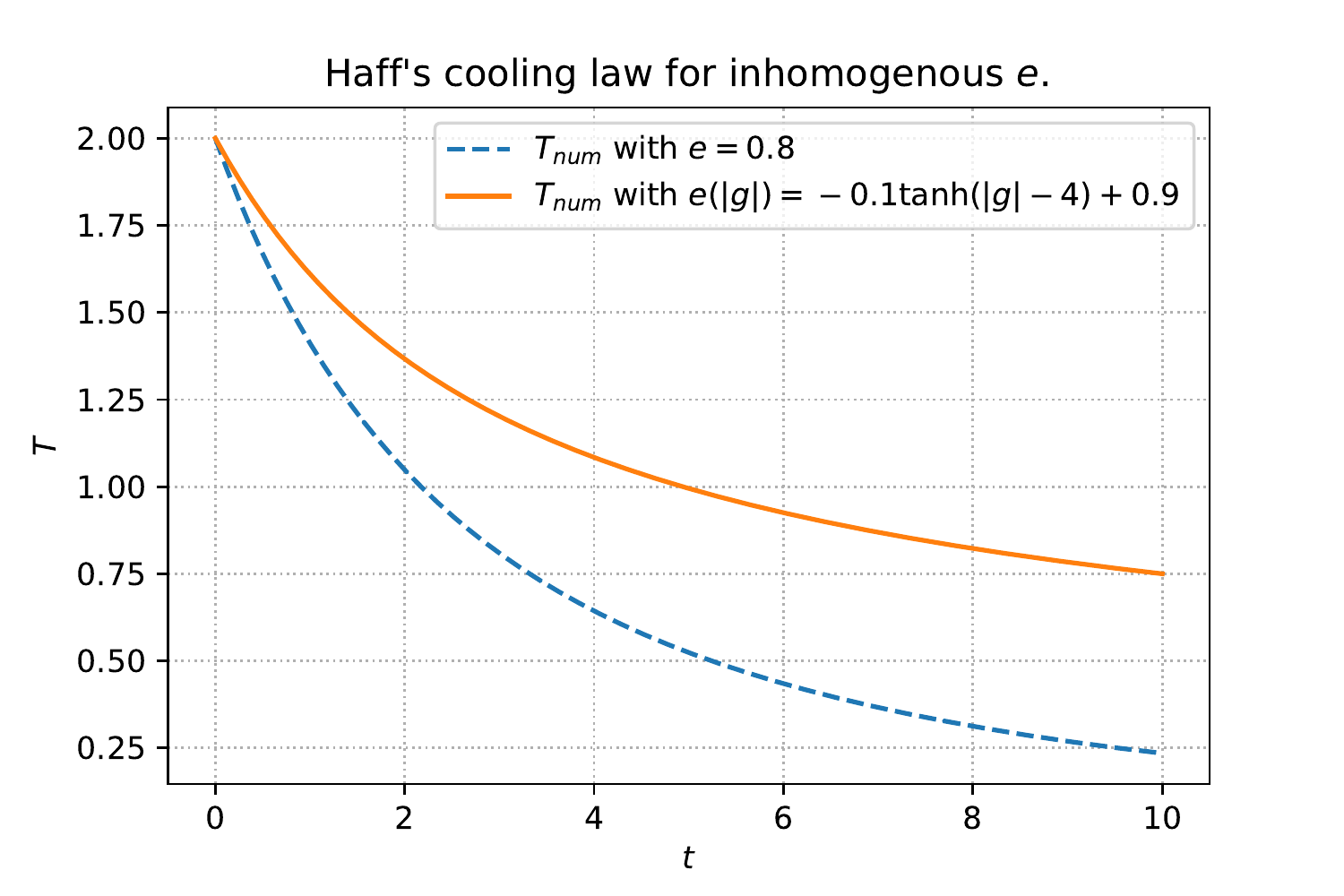}
  \caption{\textbf{Test 3.} Haff's cooling law with Maxwellian initial
    data~\eqref{eq:maxw_init_3d}. Left: plot of inhomogeneous $e$. Right: comparison of temperature between constant $e=0.8$ (dash line) and $e(|g|=\rho) = -0.1\tanh(\rho-4)+0.9$.
    Numerical parameters: $N_v^3=32\times 32\times 32$, $N_{\rho}=30$,
    $M_{\hat{g}}=32$, $R=8$, $L = 5(3+\sqrt{2})$ and $\Delta t=0.01$.}
  \label{fig:var_e_08}
\end{figure}

\begin{figure}
  \centering
  \includegraphics[width = .48\textwidth]{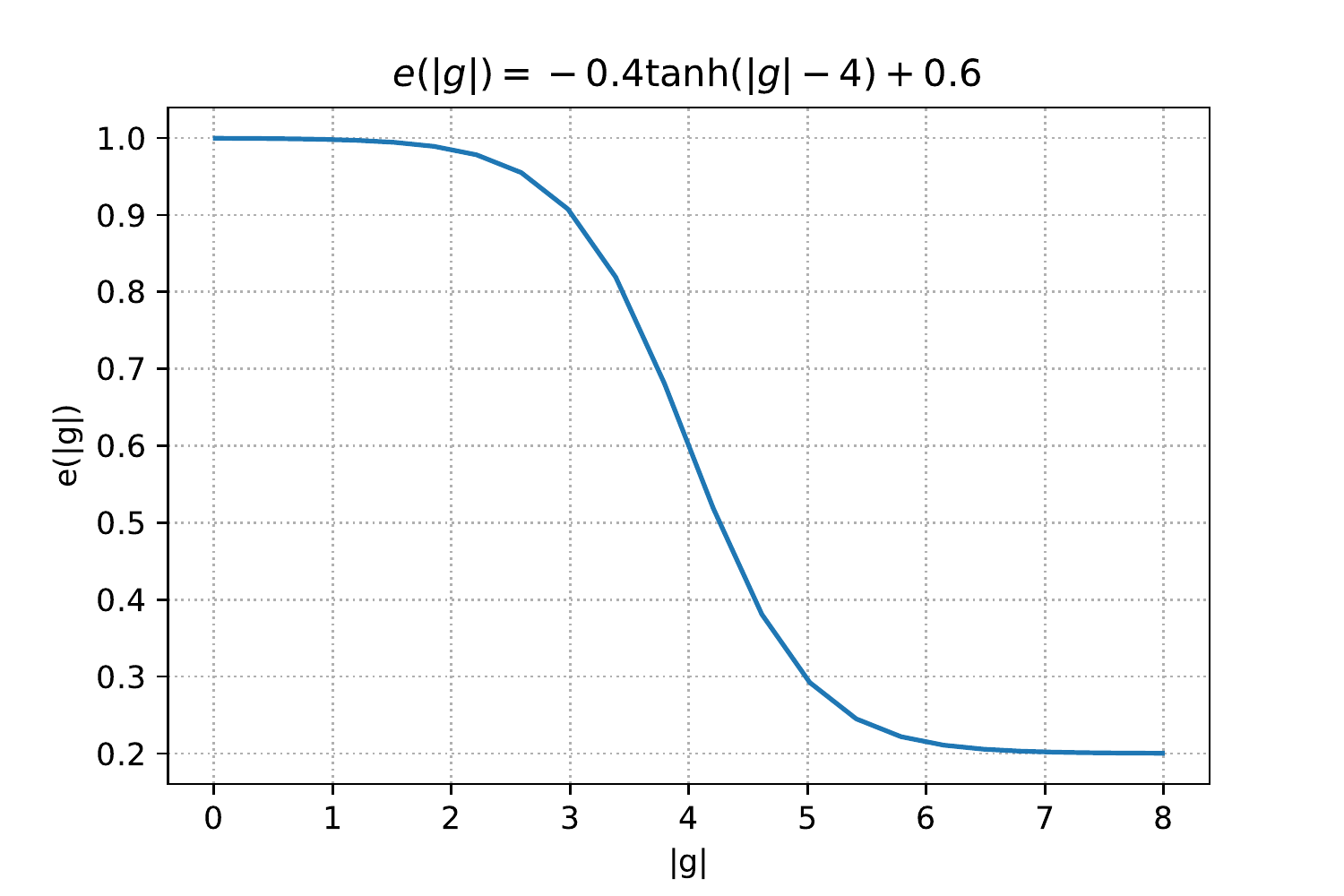}
  \includegraphics[width = .48\textwidth]{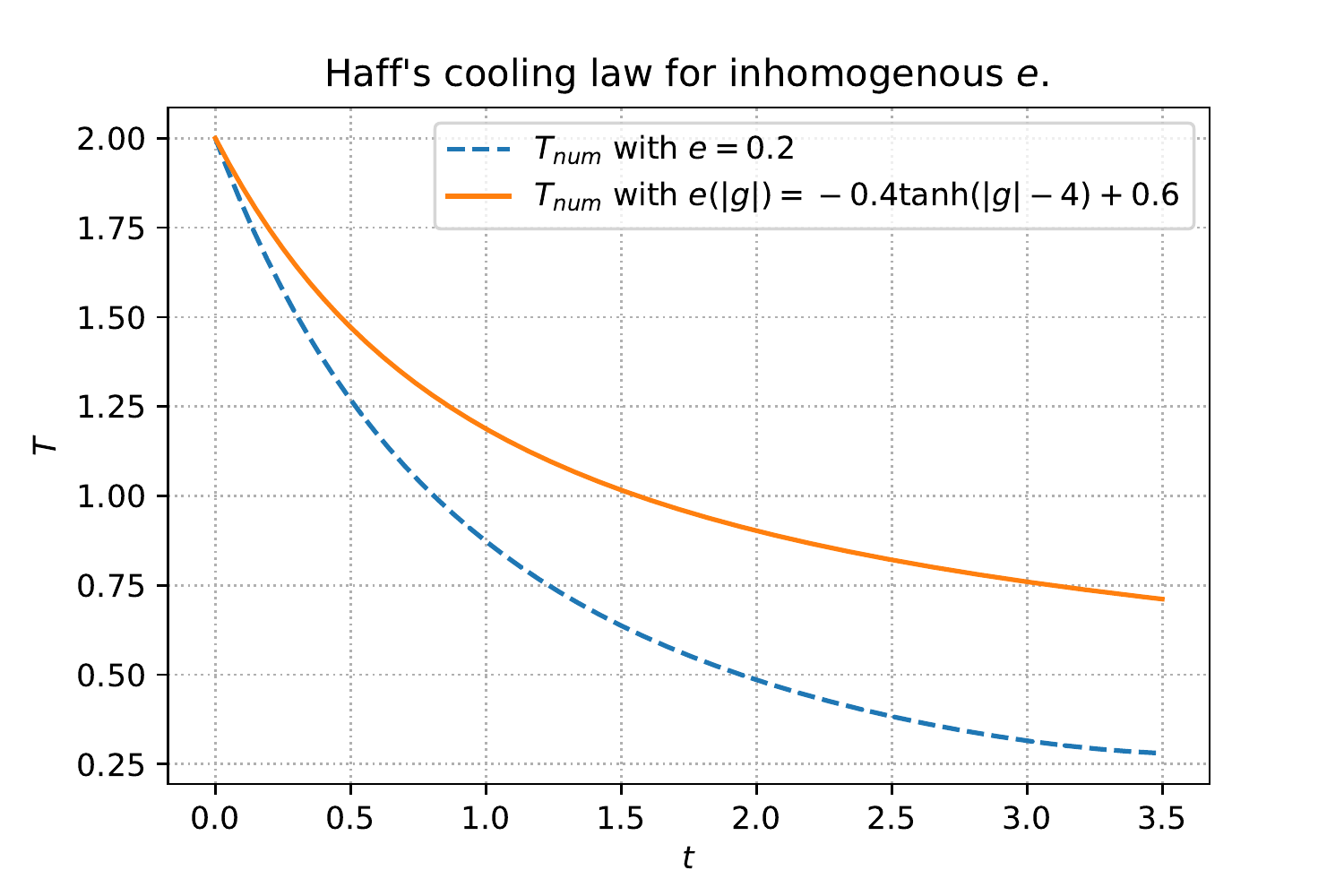}
  \caption{\textbf{Test 3.} Haff's cooling law with Maxwellian initial
    data~\eqref{eq:maxw_init_3d}. Left: plot of inhomogeneous $e$. Right: comparison of temperature between constant $e=0.2$ (dash line) and $e(|g|=\rho) = -0.4\tanh(\rho-4)+0.6$.
    Numerical parameters: $N_v^3=32\times 32\times 32$, $N_{\rho}=30$,
    $M_{\hat{g}}=32$, $R=8$, $L = 5(3+\sqrt{2})$ and $\Delta t=0.01$.}
  \label{fig:var_e_02}
\end{figure}

The last test is more related to physics in real world, by simulating the so-called
``Haff's cooling Law'' \eqref{eq:HaffLaw}. We use the following initial data which is a Maxwellian with nonzero bulk velocity:
\begin{equation}\label{eq:maxw_init_3d}
  f_0(v) = \frac{\rho_0}{(2\pi T_0)^{3/2}}e^{-(v - u_0)^2},
\end{equation}
where $\rho_0 = 1$, $T_0 = 2$ and $u_0 = (0.5, -0.5, 0)^T$. We consider the hard spheres collision kernel in $3$D, namely
\[B = \frac{1}{4\pi}|g|.\] 
In the first two tests, we consider a realistic set-up where the restitution coefficient $e$ depends on the distance of the relative velocity, i.e., $e$ is a function of $\rho=|g|$ instead of a constant,
\[
  e(\rho) = \frac{e_0 - 1}{2}\tanh(\rho-4) + \frac{e_0+1}{2},
\]
where $0<e_0<1$. This allows to mimics the physically relevant visco-elastic hard spheres case (see also \eqref{defRestiCoeffViscoEl}). We numerically evaluate the temperature and the results for $e_0=0.8$ and $e_0=0.2$ are shown in Fig.~\ref{fig:var_e_08} and Fig.~\ref{fig:var_e_02}. Compared with the cases where $e$ is constant, we observe a slight slower decay of the temperature.

\begin{figure}
  \centering
  \includegraphics[width = .5\textwidth]{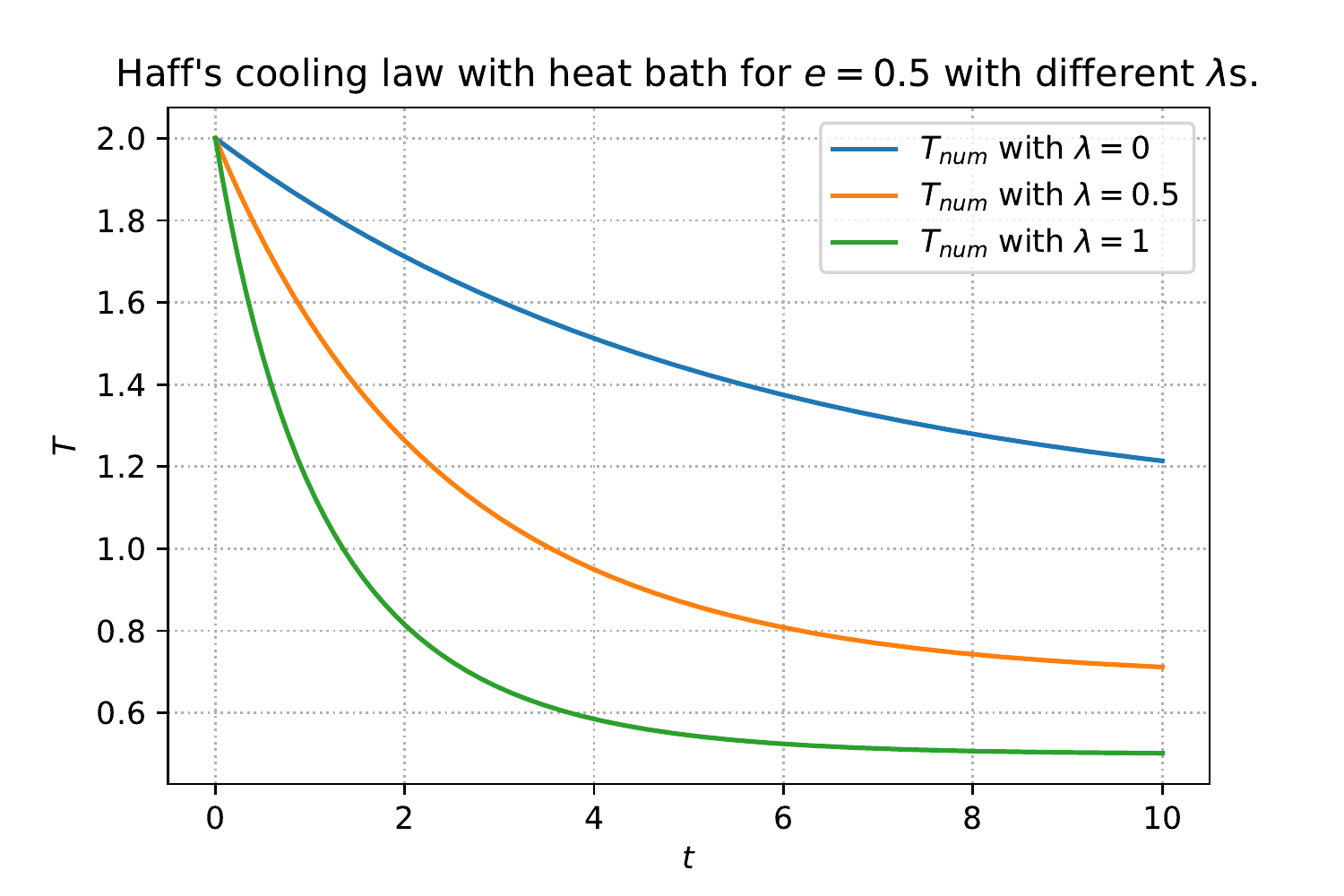}
  \caption{\textbf{Test 3.} Haff's cooling law with heat bath for different variable hard spheres exponent $\lambda$s and Maxwellian initial
    data~\eqref{eq:maxw_init_3d}. The heat bath ($\tau = 0.1$).
    Numerical parameters: $N_v^3=32\times 32\times 32$, $N_{\rho}=30$,
    $M_{\hat{g}}=32$, $R=8$, $L = 5(3+\sqrt{2})$ and $\Delta t=0.01$.}
  \label{fig:var_gamma}
\end{figure}

Another parameter that may affect the decay rate of temperature is the variable hard spheres exponent $\lambda$ from \eqref{defHSkernelAnomalous}. In Fig.~\ref{fig:var_gamma} we show that, in the presence of heat bath, for $e=0.5$ but with $\lambda=1$ (hard spheres), $\lambda=0.5$ and $\lambda=0$ (Maxwellian molecules), the decay rate of temperature will decrease after certain time (notice the slopes after $t=5$). 

\begin{figure}
  \centering
  \includegraphics[width = .48\textwidth]{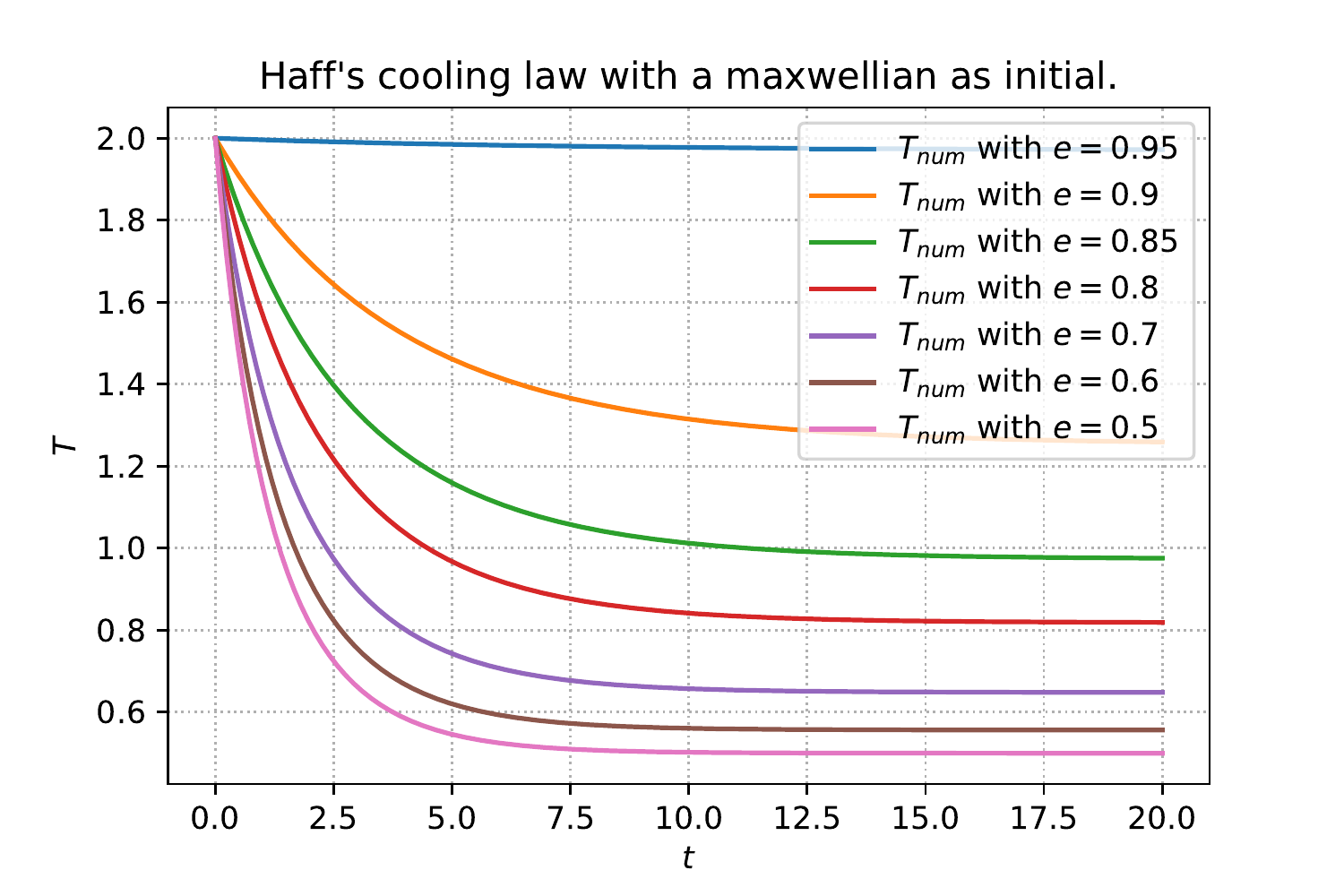}
  \includegraphics[width = .48\textwidth]{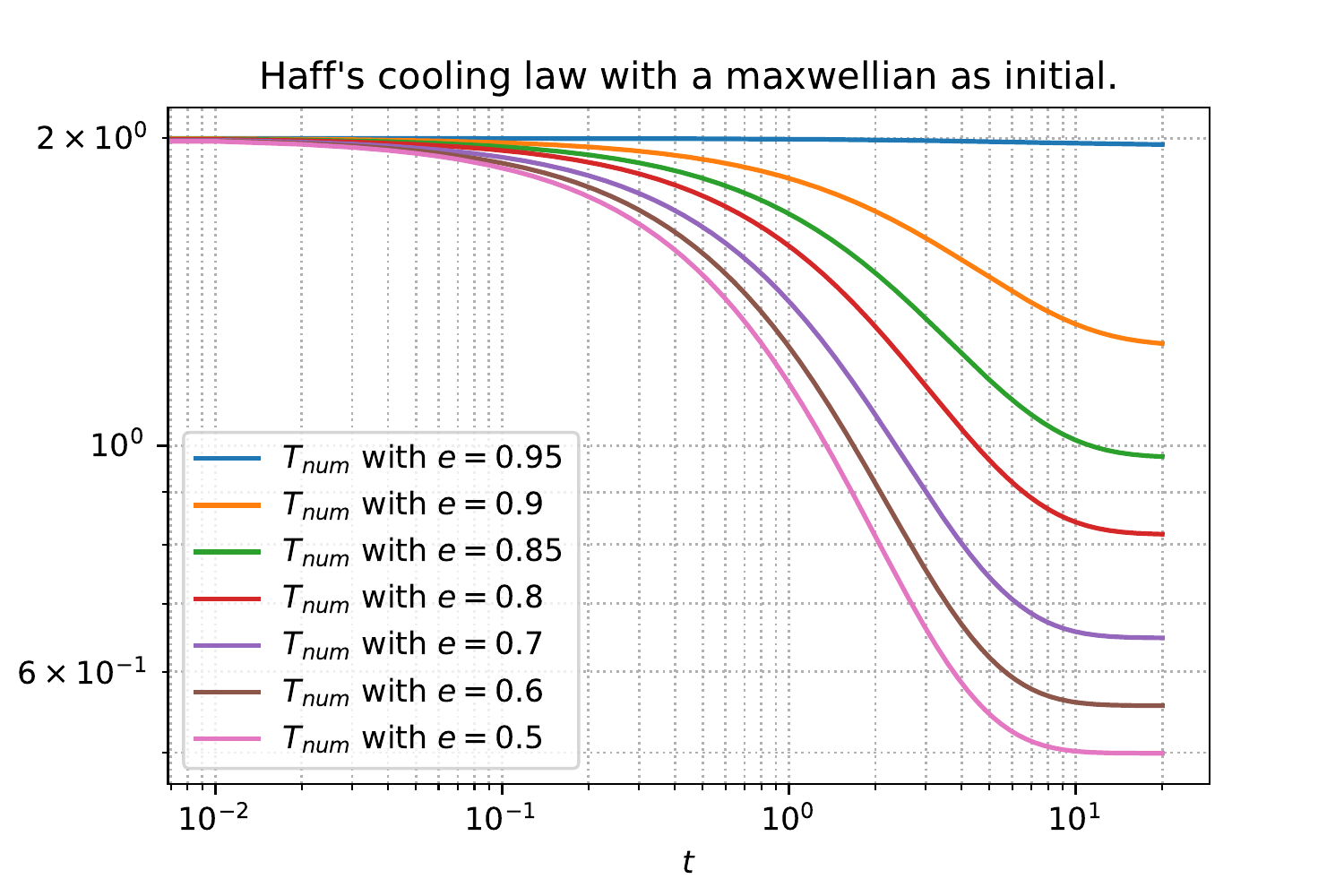}
  \caption{\textbf{Test 3.} Heated Haff's cooling law with Maxwellian initial
    data~\eqref{eq:maxw_init_3d}. Left: regular plot. Right: log-log plot.
    Numerical parameters: $N_v^3=32\times 32\times 32$, $N_{\rho}=30$,
    $M_{\hat{g}}=32$, $R=8$, $L = 5(3+\sqrt{2})$ and $\Delta t=0.01$.}
  \label{fig:heated_Haff_cooling_maxw}
\end{figure}

Finally, with the heat bath $\tau = 0.1$, we
numerically evaluate the temperature up to time
$t_\text{final}=20$ for various values of restitution coefficients.
The time evolution of $T$ is shown in Fig.~\ref{fig:heated_Haff_cooling_maxw} where one can observe the
transition of decays from $e=0.5$ to $e=0.95$ (near elastic case).

\bibliographystyle{acm}
\bibliography{biblio,mybibfile}

\end{document}